\documentclass[11pt]{article}
\usepackage{amssymb, amsmath, amsthm, mathrsfs, ascmac, color}
\usepackage[pdftex]{graphicx} 
\usepackage[subrefformat=parens]{subcaption}
\usepackage[top=3cm, bottom=3cm, left=3cm,right=3cm]{geometry}
\usepackage{comment}

\theoremstyle{plain}
\newtheorem{thm}{Theorem}[section]

\newtheorem{rmk}{Remark}
\DeclareFontFamily{U}{mathx}{}
\DeclareFontShape{U}{mathx}{m}{n}{<-> mathx10}{}
\DeclareSymbolFont{mathx}{U}{mathx}{m}{n}
\DeclareMathAccent{\widecheck}{0}{mathx}{"71}
\makeatletter

\@addtoreset{equation}{section}
\makeatother

\newcommand{\dd}{\mathrm d}
\newcommand{\pd}{\partial}
\newcommand{\ee}{\mathrm e}
\newcommand{\EE}{\mathbb E}
\newcommand{\PP}{\mathbf{P}}
\newcommand{\OO}{\mathcal O}
\newcommand{\oo}{{\scriptsize\text{$\mathcal O$}}}

\newcommand{\VV}{\mathrm{Var}}
\newcommand{\cov}{\mathrm{Cov}}

\newcommand{\pto}{\stackrel{\PP}{\to}}
\newcommand{\dto}{\stackrel{d}{\to}}
\newcommand{\TT}{\mathsf T}

\newcommand{\tand}{\widetilde \land}

\allowdisplaybreaks
\title{\textbf{Estimation for linear parabolic SPDEs in two space dimensions
with unknown damping parameters
}}
\date{}

\author{\textbf{Yozo Tonaki}\thanks{Graduate School of Engineering Science, The University of Osaka, Toyonaka, Japan}
\thanks{Center for Mathematical Modeling and Data Science (MMDS), The University of Osaka, Toyonaka, Japan} 
\footnote{email: \texttt{y.tonaki.es@osaka-u.ac.jp}}
\and \textbf{Yusuke Kaino}\thanks{Graduate School of Maritime Sciences, Kobe University}
\and \textbf{Masayuki Uchida}$^{* \dag}$\thanks{CREST, Japan Science and Technology Agency, Tokyo, Japan}
}

\begin{document}
\bibliographystyle{plain}
\maketitle

\begin{abstract}
We study parametric estimation 
for second order linear parabolic stochastic partial differential equations (SPDEs)
in two space dimensions driven by two types of $Q$-Wiener processes
based on high frequency spatio-temporal data.
First, we give estimators for damping parameters of the $Q$-Wiener processes
of the SPDE using realized quadratic variations based on temporal and spatial increments.
We next propose minimum contrast estimators of four coefficient parameters in the SPDE
and obtain estimators of the rest of unknown parameters in the SPDE
using an approximate coordinate process.
We also examine  numerical simulations of the proposed estimators.
\begin{center}
\textbf{Keywords and phrases}
\end{center}
Damping factor,
high frequency spatio-temporal data,
linear parabolic stochastic partial differential equations,
parametric estimation, 
$Q$-Wiener process.
\end{abstract}

\section{Introduction}
We consider parametric estimation for a linear parabolic 
stochastic partial differential equation (SPDE)
\begin{align}
\dd X_t(y,z)
&=\biggl\{
\theta_2
\biggl(\frac{\pd^2}{\pd y^2}+\frac{\pd^2}{\pd z^2}\biggr)
+\theta_1\frac{\pd}{\pd y} 
+\eta_1\frac{\pd}{\pd z} 
+\theta_0 
\biggr\} X_t(y,z) \dd t
\nonumber
\\
&\qquad
+\sigma \dd W_t^{Q}(y,z),
\quad (t,y,z) \in [0,1] \times D
\label{2d_spde}
\end{align}
with an initial value $X_0$
and the Dirichlet boundary condition
$X_t(y,z) = 0$, $(t,y,z) \in [0,1] \times \pd D$,
where $D = (0,1)^2$, $W_t^{Q}$ is a $Q$-Wiener process in a Sobolev space on $D$,
$\theta=(\theta_0, \theta_1, \eta_1, \theta_2)$ 
and $\sigma$ are unknown parameters, 
the parameter space of $(\theta,\sigma)$ is a compact convex subset of 
$\mathbb R^3 \times (0,\infty)^2$ and the true value 
$(\theta^*, \sigma^*) 
=(\theta_0^*, \theta_1^*, \eta_1^*, \theta_2^*, \sigma^*)$ 
belongs to the interior of the parameter space.

SPDEs, which can describe real world phenomena in more detail, 
become a popular modeling tool
and are expected to be applied in many fields,
such as geophysical fluid dynamics, mathematical finance, and environmental science.
In particular, parabolic SPDEs in two space dimensions 
can be used for modeling of sea surface temperature anomalies, 
climate energy balance and plutonium contamination, 
see Piterbarg and Ostrovskii \cite{Piterbarg_Ostrovskii1997},
North et al.\,\cite{North_etal2011}, 
Jones and Zhang \cite{Jones_Zhang1997} and Mohapl \cite{Mohapl2000}.
Along with growing interest in SPDEs, statistics for SPDEs have been developed,
see H\"ubner et al.\,\cite{Hubner_etal1993},
Hubner and Rozovskii \cite{Hubner_Rozovskii1995}, 
Piterbarg and Rozovskii \cite{Piterbarg_Rozovskii1997},
Ibragimov and Khasminskii \cite{Ibragimov_Khasminskii1999},
Prakasa Rao \cite{PrakasaRao2001},
Lototsky \cite{Lototsky2003}, 
Markussen \cite{Markussen2003},
Cialenco and Glatt-Holtz \cite{Cialenco_Glatt-Holtz2011}, 
Cialenco \cite{Cialenco2018},
Mahdi Khalil and Tudor \cite{MahdiKhalil_Tudor2019},
Assaad and Tudor \cite{Assaad_Tudor2021},
Pang and Cao \cite{Pang_Cao2021},
Tudor \cite{Tudor2022},
and Benth et al.\,\cite{Benth_etal2024}.

Statistical inference for second order parabolic SPDEs based on high frequency data
has been studied by many researchers, see for example, 
Bibinger and Trabs \cite{Bibinger_Trabs2020},
Chong \cite{Chong2020},
Cialenco and Huang \cite{Cialenco_Huang2020},
Kaino and Uchida \cite{Kaino_Uchida2020,Kaino_Uchida2021},
Hildebrandt and Trabs \cite{Hildebrandt_Trabs2021},
Bibinger and Bossert \cite{Bibinger_Bossert2023},
Gamain and Tudor \cite{Gamain_Tudor2023},
Tonaki et al.\,\cite{TKU2023, TKU2024a, TKU2024b, TKU2024arXiv1, 
TKU2024arXiv2, TKU2024arXiv3}, and Gaudlitz and Reiss \cite{Gaudlitz_Reiss2023}.
Tonaki et al.\,\cite{TKU2023,TKU2024b} studied the estimation 
of the coefficient parameters in SPDE \eqref{2d_spde} driven by 
a $Q$-Wiener process given by \eqref{QW1} or \eqref{QW2} below
with a known damping parameter $\alpha$.
Tonaki et al.\,\cite{TKU2023} set $\alpha > 0$ 
in order to guarantee that there exists a mild solution 
of the SPDE in two space dimensions with a finite $L^2$-estimate. 
Bossert \cite{Bossert2024} considered parametric estimation for 
second order parabolic SPDEs in multiple space dimensions and,
in particular, constructed an estimator of 
the damping parameter $\alpha$ of a $Q$-Wiener process based on realized volatilities.
\cite{Bossert2024} showed that the estimator of the damping parameter 
is asymptotically normal with the rate 
$\sqrt{mn}$, where $n$ and $m$ denote the number of temporal and spatial observations, 
respectively, and $m = \OO(n^\rho)$ with for some $\rho \in (0,(1-\alpha)/4)$
($\alpha \in (0,1)$) in the case of the SPDE in two space dimensions.
\cite{Bossert2024} also stated that one can construct consistent estimators 
of the curvature parameter and the normalized volatility 
using this estimator even if the damping parameter $\alpha$ is unknown.

In this paper, we propose an estimator of the damping parameter $\alpha$ 
of the $Q$-Wiener process \eqref{QW1} or \eqref{QW2} below 
based on temporal and spatial realized quadratic variations.
We then substitute the estimator of $\alpha$ into the contrast function 
proposed by Tonaki et al.\,\cite{TKU2024b} to obtain estimators 
of the coefficient parameters when $\alpha$ is unknown.
The main purpose of this paper is to show that the estimator of $\alpha$ is 
bounded in probability at the parametric rate $\sqrt{m N}$
and the estimators of the coefficient parameters are 
bounded in probability at the rate $\sqrt{m N}/\log(N)$,
where $m$ is the number of spatially thinned data, 
$N$ is the number of temporal observations, $m = \OO(N)$ and $N = \OO(m)$.

This paper is organized as follows.
Section \ref{sec2} presents some notation.
In Section \ref{sec3}, we provide an estimator of the damping parameter $\alpha$ 
of the $Q$-Wiener processes given by \eqref{QW1} and \eqref{QW2}, 
and show the asymptotic properties of the estimator of $\alpha$.
In Section \ref{sec4}, we first obtain estimators of the coefficient parameters 
$\theta_2$, $\theta_1$, $\eta_1$ and $\sigma^2$ in SPDE \eqref{2d_spde} 
in line with the approach of \cite{TKU2024b}.
We then propose estimators of the rest of unknown parameters based on 
an approximate coordinate process of SPDE \eqref{2d_spde}.
In Section \ref{sec5}, we examine numerical simulations of the proposed estimators.
Section \ref{sec6} is devoted to the proofs of our results.

\section{Preliminaries}\label{sec2}
We define an operator $A_\theta$ by
\begin{equation*}
-A_\theta = 
\theta_2\biggl(\frac{\pd^2}{\pd y^2} + \frac{\pd^2}{\pd z^2} \biggr)
+ \theta_1\frac{\pd}{\pd y} + \eta_1\frac{\pd}{\pd z} + \theta_0.
\end{equation*}
The eigenfunctions $e_{l_1,l_2}$ of the operator $A_\theta$ 
and the corresponding eigenvalues $\lambda_{l_1,l_2}$ are given by
\begin{equation*}
e_{ l_1, l_2}(y,z)=e_{ l_1}^{(1)}(y)e_{ l_2}^{(2)}(z),
\quad
\lambda_{ l_1, l_2}=\theta_2(\pi^2( l_1^2+ l_2^2)+\Gamma)
\end{equation*}
for $ l_1, l_2 \ge 1$ and $y,z\in [0,1]$, where
\begin{equation*}
e_{ l_1}^{(1)}(y)
=\sqrt 2 \sin(\pi l_1 y) \ee^{-\kappa y/2},
\quad
e_{ l_2}^{(2)}(z)
=\sqrt 2 \sin(\pi l_2 z) \ee^{-\eta z/2},
\end{equation*}
\begin{equation*}
\kappa=\frac{\theta_1}{\theta_2}, 
\quad
\eta=\frac{\eta_1}{\theta_2},
\quad
\Gamma=-\frac{\theta_0}{\theta_2} +\frac{\kappa^2+\eta^2}{4}.
\end{equation*}
We assume $\lambda_{1,1}^* > 0$ 
so that $A_\theta$ is a positive definite and self-adjoint operator.
The eigenfunctions $\{e_{l_1,l_2}\}_{l_1, l_2 \ge 1}$ 
are orthonormal with respect to the weighted $L^2$-inner product
\begin{equation*}
\langle u, v \rangle
= \int_0^1 \int_0^1 u(y,z)v(y,z)\ee^{\kappa y +\eta z} \dd y \dd z, 
\quad 
\| u \| = \sqrt{\langle u, u \rangle}
\end{equation*}
for $u, v \in L^2(D)$.

Let $\{w_{ l_1, l_2}\}_{ l_1, l_2 \ge 1}$ be 
independent $\mathbb R$-valued standard Brownian motions.
Let $\alpha \in (0,2)$ be an unknown damping parameter.
Let $[\underline \alpha, \overline \alpha] \subset (0,2)$ be
the parameter space of $\alpha$ and the true value 
$\alpha^*$ belongs to $(\underline \alpha, \overline \alpha)$.
We consider the following two types of $Q$-Wiener processes introduced 
by \cite{TKU2023}. See \cite{TKU2023} for details on these two processes.
\begin{enumerate}
\item[(a)]
\textit{$Q_1$-model.}
The $Q_1$-Wiener process is given by 
\begin{equation}\label{QW1}
W_t^{Q_1} = \sum_{ l_1, l_2\ge1} \lambda_{ l_1, l_2}^{-\alpha/2} 
e_{ l_1, l_2} w_{ l_1, l_2}(t),
\quad t \ge 0.
\end{equation}

\item[(b)]
\textit{$Q_2$-model.}
The $Q_2$-Wiener process is given by
\begin{equation}\label{QW2}
W_t^{Q_2} = \sum_{ l_1, l_2\ge1} \mu_{ l_1, l_2}^{-\alpha/2} 
e_{ l_1, l_2} w_{ l_1, l_2}(t),
\quad t \ge 0,
\end{equation}
where $\mu_{l_1,l_2} = \pi^2(l_1^2 +l_2^2) +\mu_0$,
$\mu_0 \in (-2\pi^2, \infty)$ is an unknown parameter,
the parameter space of $\mu_0$ is a compact convex subset of $(-2\pi^2, \infty)$
and the true value $\mu_0^*$ belongs to the interior of the parameter space. 
\end{enumerate}

There exists a unique mild solution of SPDE \eqref{2d_spde}, which is given by
\begin{equation*}
X_t=\ee^{-t A_\theta}X_0 +\epsilon \int_0^t \ee^{-(t-s)A_\theta}\dd W_s^{Q},
\end{equation*}
where $\ee^{-t A_\theta} u 
= \sum_{ l_1, l_2 \ge 1} \ee^{-\lambda_{ l_1, l_2}t}
\langle u, e_{ l_1, l_2}\rangle e_{ l_1, l_2}$ for $u \in L^2(D)$.
$X_t$ is then decomposed as follows.
\begin{equation*}
X_t(y,z)=\sum_{ l_1, l_2\ge1} 
x_{ l_1, l_2}(t)
e_{ l_1}^{(1)}(y)e_{ l_2}^{(2)}(z),
\quad t \ge 0, \ y,z \in [0,1], 
\end{equation*}
where $x_{ l_1, l_2}(t) = \langle X_t, e_{ l_1, l_2} \rangle$.

We assume that a mild solution is discretely observed on 
the grid $(t_i,y_j,z_k) \in [0,1]^3$ with 
\begin{equation*}
t_i = i\Delta := i/N, \quad 
y_j = j/M_1, \quad
z_k = k/M_2
\end{equation*}
for $i=0,\ldots,N$, $j=0,\ldots,M_1$ and $k=0,\ldots,M_2$.
That is, we have discrete observations 
$\mathbb X_{\mathbf{M},N} = \{X_{t_i}(y_j,z_k)\}$ for 
$i=0,\ldots,N$, $j=0,\ldots,M_1$, $k=0,\ldots, M_2$ and $\mathbf{M}=(M_1, M_2)$.
For $c \in [0,1/2)$, 
$m_1 \le M_1$, $m_2 \le M_2$, $\mathbf{m} = (m_1,m_2)$ and $n \le N$,
we will write the thinned data obtained from $\mathbb X_{\mathbf{M},N}$ as 
\begin{equation*}
\mathbb X_{\mathbf{m},n}^{(c)} 
= \{ X_{\widetilde t_i}(\widetilde y_j,\widetilde z_k) \},
\quad
i=0,\ldots,n,\ j=0,\ldots,m_1,\ k=0,\ldots, m_2
\end{equation*}
with 
\begin{equation*}
\widetilde t_i = i \Delta_n = \biggl \lfloor \frac{N}{n} \biggr \rfloor \frac{i}{N},
\quad
\widetilde y_j = c + j \delta, 
\quad
\widetilde z_k = c + k \delta,
\end{equation*}
where $m_1 = m_2$, $\delta = \frac{1-2c}{m_1}$.

For a sequence $\{a_n\}$, 
we write $a_n \equiv a$ if $a_n = a$ for some $a \in \mathbb R$ and all $n$.
For $h \in (0,1)$, $a, b, c >0$ and $d \in \mathbb R$, we define
\begin{equation*}
h^{a \tand b} =  
\begin{cases}
h^{a}, & a < b,
\\
-h^{b} \log (h), & a = b,
\\
h^{b}, & a > b,
\end{cases}
\end{equation*}
and $h^{d+ c(a \tand b)} = h^d \cdot (h^c)^{a \tand b}$.
Furthermore, for $L \in (1,\infty)$, we write
\begin{equation*}
L^{a \tand b} = \frac{1}{(1/L)^{a \tand b}} =
\begin{cases}
L^{a}, & a < b,
\\
L^{b} /\log (L), & a = b,
\\
L^{b}, & a > b.
\end{cases}
\end{equation*}

The symbols $\pto$ and $\dto$ represent convergence in probability and convergence in distribution, respectively.

\section{Estimation of the damping parameter $\alpha$}\label{sec3}
Since both models can be discussed in the same way, 
we only consider the $Q_1$-model.

Let $b \in (0,1/2)$.
Suppose that we have two thinned data $\mathbb X_{\mathbf{m},N}^{(b)}$
and $\mathbb X_{\mathbf{m}',N'}^{(b)}$ obtained from $\mathbb X_{\mathbf{M},N}$,
where $\mathbf{m}' = \mathbf{m}/2$, $N' = N/4$, 
$m = m_1 m_2 = \OO(N)$ and $N = \OO(m)$.
Note that
\begin{equation*}
\delta' = \frac{1-2b}{m_1'} = 2 \times \frac{1-2b}{m_1} = 2 \delta,
\quad
\Delta' = \Delta_{N'} = 
\biggl \lfloor \frac{N}{N'} \biggr \rfloor \frac{1}{N} 
= 4 \times \frac{1}{N} = 4 \Delta.
\end{equation*}
For the thinned data $\mathbb X_{\mathbf{m},N}^{(b)}$,
we define the triple increments
\begin{align*}
T_{i,j,k} X 
&= X_{t_{i}}(\widetilde y_{j},\widetilde z_{k})
-X_{t_{i}}(\widetilde y_{j-1},\widetilde z_{k})
-X_{t_{i}}(\widetilde y_{j},\widetilde z_{k-1})
+X_{t_{i}}(\widetilde y_{j-1},\widetilde z_{k-1})
\\
&\qquad
-X_{t_{i-1}}(\widetilde y_{j},\widetilde z_{k})
+X_{t_{i-1}}(\widetilde y_{j-1},\widetilde z_{k})
+X_{t_{i-1}}(\widetilde y_{j},\widetilde z_{k-1})
-X_{t_{i-1}}(\widetilde y_{j-1},\widetilde z_{k-1}).
\end{align*}
We make the following condition of the initial value $X_0$ of SPDE \eqref{2d_spde}.
\begin{description}
\item[\textbf{[A1]}]
The initial value $X_0 \in L^2(D)$ is deterministic and
$\| A_\theta^{(1+\overline \alpha)/2} X_0 \| < \infty$. 
\end{description}
Note that $\overline \alpha$ is the upper bound of the parameter space of $\alpha$.

For the triple increments $T_{i,j,k} X$ obtained from $\mathbb X_{\mathbf{m},N}^{(b)}$, 
we obtain for $r \equiv \delta/\sqrt{\Delta} \in (0,\infty)$,
\begin{equation*}
\frac{1}{m N \Delta^\alpha} \sum_{k=1}^{m_2} \sum_{j=1}^{m_1} \sum_{i=1}^N
(T_{i,j,k} X)^2
\pto g_{r,\alpha}(\vartheta)
\end{equation*}
under [A1], where $g_{r,\alpha}(\vartheta)$ is given by \eqref{pf1-1} below.
Let $T_{i,j,k}' X$ be the triple increments obtained from $\mathbb X_{\mathbf{m}',N'}^{(b)}$.
Similarly, we have for $r' \equiv \delta'/\sqrt{\Delta'} \in (0,\infty)$,
\begin{equation*}
\frac{1}{m' N' (\Delta')^\alpha} \sum_{k=1}^{m_2'} \sum_{j=1}^{m_1'} \sum_{i=1}^{N'} 
(T_{i,j,k}' X)^2
\pto g_{r',\alpha}(\vartheta).
\end{equation*}
Since $\delta' = 2 \delta$ and $\Delta' = 4 \Delta$, we obtain $r' = r$ and
\begin{equation*}
\frac{\displaystyle 
\frac{1}{m' N'  (\Delta')^\alpha} \sum_{k=1}^{m_2'} \sum_{j=1}^{m_1'} \sum_{i=1}^{N'} 
(T_{i,j,k}' X)^2}
{\displaystyle
\frac{1}{m N \Delta^\alpha} \sum_{k=1}^{m_2} \sum_{j=1}^{m_1} \sum_{i=1}^N 
(T_{i,j,k} X)^2}
\pto 1.
\end{equation*}
Since it follows that
\begin{align*}
\frac{\displaystyle\frac{1}{m' N' (\Delta')^\alpha} \sum_{k=1}^{m'_2} \sum_{j=1}^{m'_1} 
\sum_{i=1}^{N'} (T_{i,j,k}'X)^2}
{\displaystyle\frac{1}{m N \Delta^\alpha} \sum_{k=1}^{m_2} \sum_{j=1}^{m_1} \sum_{i=1}^N 
(T_{i,j,k}X)^2}
&= \frac{\Delta^\alpha}{(\Delta')^\alpha}
\times
\frac{\displaystyle\frac{1}{m' N'} \sum_{k=1}^{m'_2} \sum_{j=1}^{m'_1} 
\sum_{i=1}^{N'} (T_{i,j,k}'X)^2}
{\displaystyle\frac{1}{m N} \sum_{k=1}^{m_2} \sum_{j=1}^{m_1} \sum_{i=1}^N 
(T_{i,j,k}X)^2}
\\
&= \frac{1}{4^\alpha}
\times
\frac{\displaystyle\frac{1}{m' N'} \sum_{k=1}^{m'_2} \sum_{j=1}^{m'_1} 
\sum_{i=1}^{N'} (T_{i,j,k}'X)^2}
{\displaystyle\frac{1}{m N} \sum_{k=1}^{m_2} \sum_{j=1}^{m_1} \sum_{i=1}^N 
(T_{i,j,k}X)^2},
\end{align*}
we have
\begin{equation*}
\frac{\displaystyle\frac{1}{m' N'} \sum_{k=1}^{m'_2} \sum_{j=1}^{m'_1} 
\sum_{i=1}^{N'} (T_{i,j,k}'X)^2}
{\displaystyle\frac{1}{m N} \sum_{k=1}^{m_2} \sum_{j=1}^{m_1} \sum_{i=1}^N 
(T_{i,j,k}X)^2}
\pto 4^\alpha.
\end{equation*}
We thus define the estimator of the damping parameter $\alpha$ as follows.
\begin{equation}\label{est_alpha}
\widehat \alpha =
\log \left( \frac{\displaystyle\frac{1}{m' N'} \sum_{k=1}^{m'_2} \sum_{j=1}^{m'_1} 
\sum_{i=1}^{N'} (T_{i,j,k}' X)^2}
{\displaystyle\frac{1}{m N} \sum_{k=1}^{m_2} \sum_{j=1}^{m_1} \sum_{i=1}^N 
(T_{i,j,k}X)^2}
\right)/\log(4).
\end{equation}

We then obtain the following result.
\begin{thm}\label{th3-1}
Under [A1], as $m \to \infty$ and $N \to \infty$,
\begin{equation*}
\sqrt{m N} (\widehat \alpha -\alpha^*) = \OO_\PP(1).
\end{equation*}
\end{thm}

\begin{rmk}
\ 
\begin{itemize}
\item[(i)]
This result also holds true for the $Q_2$-model 
without changing the form of the estimator $\widehat \alpha$.

\item[(ii)]
In the above, we chose the thinned data such that 
$\delta' = 2 \delta$ and $\Delta' = 4 \Delta$,
but an estimator of $\alpha$ can be obtained
even if we choose $\delta' = p \delta$ and $\Delta' = p^2 \Delta$ for $p = 3,4,\ldots$. 
In this case, the estimator of $\alpha$ 
is given by replacing $\log(4)$ in \eqref{est_alpha} with $\log(p^2)$.

\item[(iii)]
Since we use temporal and spatial realized quadratic variations to estimate $\alpha$, 
we can estimate a wider range of $\alpha$ than Bossert \cite{Bossert2024}, 
who used temporal realized quadratic variations to estimate $\alpha$.
\end{itemize}
\end{rmk}

\section{Estimation of the coefficient parameters}\label{sec4}
\subsection{Estimation for the $Q_1$-model}\label{sec4-1}

We consider minimum contrast estimators of the coefficient parameters 
in SPDE \eqref{2d_spde} with an unknown damping parameter $\alpha$ 
based on the methodology proposed by \cite{TKU2024b}.

Let $b \in (0,1/2)$.
Suppose that we have thinned data $\mathbb X_{\mathbf{m},N}^{(b)}$
obtained from $\mathbb X_{\mathbf{M},N}$, 
where $m = m_1 m_2$, $N = \OO(m)$ and $m = \OO(N)$.
We write 
\begin{equation*}
\overline y_j = \frac{\widetilde y_{j-1}+\widetilde y_j}{2},
\quad
\overline z_k = \frac{\widetilde z_{k-1}+\widetilde z_k}{2},
\quad
j=1,\ldots, m_1, k=1,\ldots, m_2.
\end{equation*}
Let $J_0$ be the Bessel function of the first kind of order $0$.
For $r, \alpha >0$, we set
\begin{equation}\label{psi}
\psi_{r,\alpha}(\theta_2)
=\frac{2}{\theta_2 \pi}
\int_0^\infty 
\frac{1-\ee^{-x^2}}{x^{1+2\alpha}}
\biggl(
J_0\Bigl(\frac{\sqrt{2}r x}{\sqrt{\theta_2}}\Bigr)
-2J_0\Bigl(\frac{r x}{\sqrt{\theta_2}}\Bigr)+1
\biggr) \dd x.
\end{equation}
Let $\widetilde T_{i,j,k} X = T_{i,j,k}X +T_{i+1,j,k}X$
and $\vartheta = (\kappa,\eta,\theta_2,\sigma^2)$.
We define 
\begin{align*}
\mathcal K_{\mathbf{m},N} (\vartheta;\alpha) &= 
\frac{1}{m}\sum_{k=1}^{m_2}\sum_{j=1}^{m_1} 
\Biggl\{
\frac{1}{N \Delta^{\alpha}}\sum_{i=1}^{N} (T_{i,j,k}X)^2
-f_{r,\alpha}(\overline y_j, \overline z_k;\vartheta)
\Biggr\}^2
\\
&\qquad+ 
\frac{1}{m} \sum_{k=1}^{m_2} \sum_{j=1}^{m_1} 
\Biggl\{
\frac{1}{N (2\Delta)^{\alpha}}\sum_{i=1}^{N-1} (\widetilde T_{i,j,k}X)^2
-f_{r/\sqrt{2},\alpha}(\overline y_j, \overline z_k;\vartheta)
\Biggr\}^2,
\end{align*}
where $f_{r,\alpha}(y,z;\vartheta) 
= \sigma^2 \ee^{-\kappa y -\eta z} \psi_{r,\alpha}(\theta_2)$.

Let $\Xi$ be the parameter space of $\vartheta$ and a compact convex subset of 
$\mathbb R^2 \times (-\frac{r^2}{8 \log(\sqrt{2}-1)},\infty) \times(0,\infty)$.
We assume that the true value 
$\vartheta^* = (\kappa^*,\eta^*,\theta_2^*,(\sigma^*)^2)$ 
belongs to $\mathrm{Int}(\Xi)$.
For the estimator $\widehat \alpha$ from \eqref{est_alpha}, 
we define the estimator of $\vartheta$ by
\begin{equation*}
\widehat \vartheta = \underset{\vartheta \in \Xi}{\mathrm{argmin}}\,
\mathcal K_{\mathbf{m},N}(\vartheta;\widehat \alpha).
\end{equation*}

We then obtain the following result.
\begin{thm}\label{th4-1}
Under [A1], as $m \to \infty$ and $N \to \infty$,
\begin{equation*}
\frac{\sqrt{m N}}{\log(N)}(\widehat \vartheta -\vartheta^*) = \OO_\PP(1).
\end{equation*}
\end{thm}

\begin{rmk}
\ 
\begin{enumerate}
\item[(i)]
For the lower bound
\begin{equation*}
\underline{\theta_2}(r,\alpha) = 
\begin{cases}
-\frac{r^2}{8 \log(2^{\alpha/2}-1)}, & \alpha \in (0,1),
\\
0, & \alpha \in [1,2)
\end{cases}
\end{equation*}
that guarantees the identifiability of $\theta_2$ given by \cite{TKU2024b}, 
we have 
\begin{equation*}
\sup_{\alpha \in (0,2)} \underline{\theta_2}(r,\alpha)
= -\frac{r^2}{8\log(\sqrt{2}-1)}.
\end{equation*}
That is, we set the parameter space $\Xi$ 
so that the parameter $\vartheta$ is identifiable for any $\alpha \in (0,2)$.

\item[(ii)]
For the estimator $\widehat \vartheta = 
(\widehat \kappa, \widehat \eta, \widehat \theta_2, \widehat \sigma^2)$,
we can get the estimators of the coefficient parameters $\theta_1$ and $\eta_1$  
by $\widehat \theta_1 = \widehat \kappa \widehat \theta_2$ and
$\widehat \eta_1 = \widehat \eta \widehat \theta_2$, respectively.
Therefore, the estimator
$(\widehat \theta_1, \widehat \eta_1, \widehat \theta_2, \widehat \sigma^2)$
is bounded in probability at the rate $\sqrt{m N}/\log(N)$.

\item[(iii)]
Theorem \ref{th4-1} shows that the convergence rates of the estimators for the coefficient parameters 
$\theta_1$, $\eta_1$, $\theta_2$ and $\sigma^2$ 
are reduced by  $1/\log(N)$ when $\alpha$ is unknown, compared to when $\alpha$ is known
(Theorem 2.2 in \cite{TKU2024b}).
\end{enumerate}
\end{rmk}

We next study the estimator with asymptotic normality at the rate $\sqrt{n}$
using the thinned data $\mathbb X_{\mathbf{M},n}^{(0)}$, 
where $\mathbf{M} = (M_1, M_2)$ and $n \le N$.
The temporal thinning is essential to make the number of observations in space 
relatively larger than that in time, 
see Theorem \ref{th4-4} below for the balance condition between $n$ and $\mathbf{M}$.
In the $Q_1$-model, the coordinate process is given by 
\begin{equation*}
x_{ l_1, l_2}(t) = \langle X_t, e_{ l_1, l_2} \rangle =
\ee^{-\lambda_{l_1, l_2} t} \langle X_0, e_{ l_1, l_2} \rangle 
+\sigma \lambda_{l_1, l_2}^{-\alpha/2} 
\int_0^t \ee^{-\lambda_{ l_1, l_2}(t-s)} \dd w_{ l_1, l_2}(s)
\end{equation*}
and satisfies
\begin{equation}\label{OU1}
\dd x_{ l_1, l_2}(t) =
-\lambda_{ l_1, l_2} x_{ l_1, l_2}(t)\dd t
+\sigma \lambda_{ l_1, l_2}^{-\alpha/2} \dd w_{ l_1, l_2}(t),
\quad
x_{ l_1, l_2}(0) = \langle X_0, e_{ l_1, l_2} \rangle.
\end{equation}
By the definition of the inner product $\langle \cdot, \cdot \rangle$, we have 
\begin{equation*}
x_{l_1, l_2}(t)
=2\int_0^1 \int_0^1 X_t(y,z)\sin(\pi  l_1 y) \sin(\pi  l_2 z)
\ee^{(\kappa y +\eta z)/2} \dd y \dd z.
\end{equation*}
Using the operator $\Psi_{\mathbf{M}}$ defined by
$\Psi_{\mathbf{M}} f(y,z) = f(y_{j-1}, z_{k-1})$
for $(y,z) \in [y_{j-1}, y_j) \times [z_{k-1}, z_k)$, 
$j=1,\ldots, M_1$, $k=1,\ldots, M_2$ and $f \in C(D)$, 
we approximate the coordinate process $x_{l_1,l_2}(t)$ as follows.
\begin{align*}
\widehat x_{l_1, l_2}(t)
&= 2\int_0^1 \int_0^1 \Psi_{\mathbf{M}} X_t(y,z)\sin(\pi  l_1 y) \sin(\pi  l_2 z)
\ee^{(\widehat \kappa y +\widehat \eta z)/2} \dd y \dd z
\\
&=
\sum_{j=1}^{M_1} \sum_{k=1}^{M_2} X_t(y_{j-1}, z_{k-1})
\delta_j^{[y]} g_{l_1}(\widehat \kappa) \delta_k^{[z]} g_{l_2}(\widehat \eta),
\end{align*}
where
\begin{align*}
g_l(x:a) &= 
\frac{\sqrt{2}\ee^{ax/2}}{(a/2)^2+(\pi l)^2}
\biggl( \frac{a}{2} \sin(\pi l x) - \pi l \cos(\pi l x) \biggr),
\quad a, x \in \mathbb R, \ l \in \mathbb N,
\\
\delta_j^{[y]} g_l(a) &= g_l(y_{j}:a) - g_l(y_{j-1}:a),
\quad
\delta_k^{[z]} g_l(a) = g_l(z_{k}:a) - g_l(z_{k-1}:a).
\end{align*}
We define the estimator of 
$\sigma_{l_1,l_2}:=\sigma\lambda_{l_1,l_2}^{-\alpha/2}$ in \eqref{OU1} as 
\begin{equation*}
\widetilde \sigma_{l_1,l_2}^2 = 
\sum_{i=1}^n \bigl(\widetilde x_{l_1,l_2}(\widetilde t_i)
-\widetilde x_{l_1,l_2}(\widetilde t_{i-1}) \bigr)^2.
\end{equation*}
We propose the estimators of $\lambda_{l_1,l_2}$ and
$(\theta_0, \theta_1, \eta_1, \theta_2, \sigma^2)$ as follows.
\begin{equation*}
\widetilde \lambda_{l_1,l_2} 
= \biggl(\frac{\widehat \sigma^2}{\widetilde \sigma_{l_1,l_2}^2} 
\biggr)^{1/\widehat \alpha},
\quad
\widetilde \theta_0 = -\widetilde \lambda_{1,1}+
\biggl(\frac{\widehat \kappa^2 +\widehat \eta^2}{4} +2\pi^2 \biggr) \widehat \theta_2,
\end{equation*}
\begin{equation*}
\widetilde \theta_2 = \frac{\widetilde \lambda_{1,2} -\widetilde \lambda_{1,1}}{3\pi^2},
\quad
\widetilde \theta_1 = \widehat \kappa \widetilde \theta_2,
\quad
\widetilde \eta_1 = \widehat \eta \widetilde \theta_2,
\quad
\widetilde \sigma^2 = \frac{\widehat \sigma^2}{\widehat \theta_2} \widetilde \theta_2.
\end{equation*}
Let
\begin{equation*}
\mathcal J = 2
\begin{pmatrix}
J_{1,1} & J_{1,2}
\\
J_{1,2}^\TT & J_{2,2}
\end{pmatrix},
\end{equation*}
where
\begin{equation*}
J_{1,1} = \frac{(\lambda_{1,1}^*)^2}{(\alpha^*)^2}, 
\quad
J_{1,2} =\frac{(\lambda_{1,1}^*)^2}{3\pi^2\theta_2^* (\alpha^*)^2} 
(\theta_1^*, \eta_1^*, \theta_2^*, (\sigma^*)^2),
\end{equation*}
\begin{equation*}
J_{2,2} = 
\frac{2((\lambda_{1,1}^*)^2 + (\lambda_{1,2}^*)^2)}{9\pi^4 (\theta_2^*)^2 (\alpha^*)^2}
\begin{pmatrix}
(\theta_1^*)^2 & \theta_1^* \eta_1^* & \theta_1^* \theta_2^* & \theta_1^* (\sigma^*)^2
\\
\theta_1^* \eta_1^* & (\eta_1^*)^2 & \eta_1^* \theta_2^* & \eta_1^* (\sigma^*)^2
\\
\theta_1^* \theta_2^* & \eta_1^* \theta_2^* & (\theta_2^*)^2 & \theta_2^* (\sigma^*)^2
\\
\theta_1^* (\sigma^*)^2 & \eta_1^* (\sigma^*)^2 & \theta_2^* (\sigma^*)^2 & (\sigma^*)^4
\\
\end{pmatrix}
\end{equation*}
and $\TT$ denotes the transpose.
\begin{thm}\label{th4-2}
Under [A1], it follows that
\begin{itemize}
\item[(i)]
if $\frac{n}{(M_1 \land M_2)^{2(\alpha^* \tand 1)}} \to 0$, then
$\widetilde \theta_0 \pto \theta_0^*$ as 
$M_1 \land M_2 \to \infty$ and $n \to \infty$,

\item[(ii)]
if $\frac{n}{(M_1 \land M_2)^{\alpha^* \tand 1}} \to 0$, then
as $M_1 \land M_2 \to \infty$ and $n \to \infty$,
\begin{equation}\label{eq-th2}
\sqrt{n}
\begin{pmatrix}
\widetilde \theta_0 - \theta_0^*
\\
\widetilde \theta_1 - \theta_1^*
\\
\widetilde \eta_1 - \eta_1^*
\\
\widetilde \theta_2 - \theta_2^*
\\
\widetilde \sigma^2 - (\sigma^*)^2
\end{pmatrix}
\dto N (0, \mathcal J).
\end{equation}
\end{itemize}
\end{thm}

\begin{rmk}
\ 
\begin{itemize}
\item[(i)]
Since the true value $\alpha^*$ is unknown, 
we have no way to check the sufficient conditions of 
(i) and (ii) in Theorem \ref{th4-2}. 
Therefore, in practice, 
we can only check the validity of the condition 
where $\alpha^*$ is replaced by $\widehat \alpha$.

\item[(ii)]
Even when $\alpha$ is unknown, the estimators of the coefficient parameters 
based on the approximate coordinate process attain 
the same convergence rates as in the case
when $\alpha$ is known (Theorem 2.3 in \cite{TKU2024b}). 
\end{itemize}
\end{rmk}

\subsection{Estimation for the $Q_2$-model}\label{sec4-2}
We here consider SPDE \eqref{2d_spde} driven by the $Q_2$-Wiener process given by
\eqref{QW2}.
For $r, \alpha >0$, we set
$\widetilde \psi_{r,\alpha}(\theta_2) = \theta_2^\alpha \psi_{r,\alpha}(\theta_2)$ 
for $\psi_{r,\alpha}$ given by \eqref{psi}.
We define 
\begin{align*}
\widetilde {\mathcal K}_{\mathbf{m},N} (\vartheta;\alpha) &= 
\frac{1}{m}\sum_{k=1}^{m_2}\sum_{j=1}^{m_1} 
\Biggl\{
\frac{1}{N \Delta^{\alpha}}\sum_{i=1}^{N} (T_{i,j,k}X)^2
-\widetilde f_{r,\alpha}(\overline y_j, \overline z_k;\vartheta)
\Biggr\}^2
\\
&\qquad+ 
\frac{1}{m} \sum_{k=1}^{m_2} \sum_{j=1}^{m_1} 
\Biggl\{
\frac{1}{N (2\Delta)^{\alpha}}\sum_{i=1}^{N-1} (\widetilde T_{i,j,k}X)^2
-\widetilde f_{r/\sqrt{2},\alpha}(\overline y_j, \overline z_k;\vartheta)
\Biggr\}^2,
\end{align*}
where $\widetilde f_{r,\alpha}(y,z;\vartheta) 
= \sigma^2 \ee^{-\kappa y -\eta z} \widetilde \psi_{r,\alpha}(\theta_2)$.
We define the estimator of $\vartheta$ by
\begin{equation*}
\widecheck \vartheta = \underset{\vartheta \in \Xi}{\mathrm{argmin}}\,
\widetilde {\mathcal K}_{\mathbf{m},N}(\vartheta;\widehat \alpha).
\end{equation*}

We obtain the following result.
\begin{thm}\label{th4-3}
Under [A1], as $m \to \infty$ and $N \to \infty$,
\begin{equation*}
\frac{\sqrt{m N}}{\log(N)}(\widecheck \vartheta -\vartheta^*) = \OO_\PP(1).
\end{equation*}
\end{thm}

Since in the $Q_2$-model, the coordinate process satisfies
\begin{equation}\label{OU2}
\dd x_{ l_1, l_2}(t) =
-\lambda_{ l_1, l_2} x_{ l_1, l_2}(t)\dd t
+\sigma \mu_{ l_1, l_2}^{-\alpha/2} \dd w_{ l_1, l_2}(t),
\quad
x_{ l_1, l_2}(0) = \langle X_0, e_{ l_1, l_2} \rangle,
\end{equation}
we define the approximate coordinate process 
\begin{equation*}
\overline x_{l_1, l_2}(t)
= \sum_{j=1}^{M_1} \sum_{k=1}^{M_2} X_t(y_{j-1}, z_{k-1})
\delta_j^{[y]} g_{l_1}(\widecheck \kappa) \delta_k^{[z]} g_{l_2}(\widecheck \eta).
\end{equation*}
Let $\varsigma_{l_1,l_2}=\sigma \mu_{l_1,l_2}^{-\alpha/2}$. 
The estimator of $\varsigma_{l_1,l_2}^2$ is given by
\begin{equation*}
\overline \varsigma_{l_1,l_2}^2 = 
\sum_{i=1}^n \bigl(\overline x_{l_1,l_2}(\widetilde t_i)
-\overline x_{l_1,l_2}(\widetilde t_{i-1}) \bigr)^2.
\end{equation*}
We propose the estimators of $\mu_0$ and
$(\theta_1, \eta_1, \theta_2, \sigma^2)$ as follows.
\begin{equation*}
\overline \sigma^2 = 
\Biggl\{
3\pi^2 \biggl( \frac{1}{\overline \varsigma_{1,2}^{2/\widehat \alpha}}
-\frac{1}{\overline \varsigma_{1,1}^{2/\widehat \alpha}}
\biggr)^{-1}
\Biggr\}^{\widehat \alpha},
\quad
\overline \mu_0 
= \biggl(\frac{\widecheck \sigma^2}{\overline \varsigma_{1,1}^2}
\biggr)^{1/\widehat \alpha} -2 \pi^2,
\end{equation*}
\begin{equation*}
\overline \theta_2 = \frac{\widecheck \theta_2}{\widecheck \sigma^2}\overline \sigma^2,
\quad
\overline \theta_1 = \widecheck \kappa \overline \theta_2,
\quad
\overline \eta_1 = \widecheck \eta \overline \theta_2,
\quad
\end{equation*}
Let
\begin{equation*}
\mathcal K = 2
\begin{pmatrix}
K_{1,1} & K_{1,2}
\\
K_{1,2}^\TT & K_{2,2}
\end{pmatrix},
\end{equation*}
where
\begin{equation*}
K_{1,1} = \frac{(\mu_{1,1}^*)^2}{(\alpha^*)^2}, 
\quad
K_{1,2} =\frac{(\mu_{1,1}^*)^2}{3\pi^2\theta_2^* (\alpha^*)^2} 
(\theta_1^*, \eta_1^*, \theta_2^*, (\sigma^*)^2),
\end{equation*}
\begin{equation*}
K_{2,2} = 
\frac{2((\mu_{1,1}^*)^2 + (\mu_{1,2}^*)^2)}{9\pi^4}
\begin{pmatrix}
(\theta_1^*)^2 & \theta_1^* \eta_1^* & \theta_1^* \theta_2^* & \theta_1^* (\sigma^*)^2
\\
\theta_1^* \eta_1^* & (\eta_1^*)^2 & \eta_1^* \theta_2^* & \eta_1^* (\sigma^*)^2
\\
\theta_1^* \theta_2^* & \eta_1^* \theta_2^* & (\theta_2^*)^2 & \theta_2^* (\sigma^*)^2
\\
\theta_1^* (\sigma^*)^2 & \eta_1^* (\sigma^*)^2 & \theta_2^* (\sigma^*)^2 & (\sigma^*)^4
\\
\end{pmatrix}.
\end{equation*}
\begin{thm}\label{th4-4}
Under [A1], it follows that
\begin{itemize}
\item[(i)]
if $\frac{n}{(M_1 \land M_2)^{2(\alpha^* \tand 1)}} \to 0$, then
$\overline \mu_0 \pto \mu_0^*$ as 
$M_1 \land M_2 \to \infty$ and $n \to \infty$.

\item[(ii)]
if $\frac{n^2}{(M_1 \land M_2)^{2(\alpha^* \tand 1)}} \to 0$, then
as $M_1 \land M_2 \to \infty$ and $n \to \infty$,
\begin{equation*}
\sqrt{n}
\begin{pmatrix}
\overline \mu_0 - \mu_0^*
\\
\overline \theta_1 - \theta_1^*
\\
\overline \eta_1 - \eta_1^*
\\
\overline \theta_2 - \theta_2^*
\\
\overline \sigma^2 - (\sigma^*)^2
\end{pmatrix}
\dto N (0, \mathcal K).
\end{equation*}
\end{itemize}
\end{thm}

\section{Simulation results}\label{sec5}
We generated the numerical solution of SPDE \eqref{2d_spde} 
with the $Q_1$-Wiener process
by using the truncation method with the Euler-Maruyama scheme: 
\begin{equation*}
\widetilde X_{t_{i}}(y_j, z_k)
= \sum_{l_1=1}^{L_1} \sum_{l_2=1}^{L_2} 
x_{l_1,l_2}(t_{i}) e_{l_1,l_2}(y_j, z_k), 
\quad i = 1,\ldots, N, j = 1,\ldots, M_1, k = 1,\ldots, M_2.
\end{equation*}
In this simulation, the true values of parameters
$(\theta_0^*, \theta_1^*,\eta_1^*, \theta_2^*, \sigma^*) = (0,0.2,0.2,0.2,1)$.
We set that $N = 10^3$, $M_1 = M_2 = 200$, $L_1 = L_2 = 10^4$,
$X_0 = 0$ and $\alpha^* = 0.5$.
The number of Monte Carlo iterations is 250.

\begin{table}[h]
\caption{Simulation results of $\widehat{\alpha}$, $\widehat{\theta}_0$, 
$\widehat{\theta}_1$, $\widehat{\eta}_1$, $\widehat{\theta}_2$ 
and $\widehat{\sigma}^2$.  \label{tab1}}
\begin{center}
\begin{tabular}{cc|cccccc} \hline
		& &$\widehat{\alpha}$ & $\widehat{\theta}_0$ & $\widehat{\theta}_1$ & $\widehat{\eta}_1$ & $\widehat{\theta}_2$ & $\widehat{\sigma}^2$
\\ \hline
&true value &0.5 &0 & 0.2& 0.2 & 0.2 & 1
\\ \hline
Case 1&mean & 0.497 &-0.368& 0.202 & 0.202 & 0.205 & 0.998
 \\
&s.d. & (0.001)&(0.972) & (0.002) & (0.002) & (0.002) & (0.007)
 \\   \hline
Case 2&mean & 0.482 &-0.275& 0.220 & 0.220 & 0.223 & 0.990
 \\
&s.d. & (0.002)&(1.088) & (0.004) & (0.005) & (0.004) & (0.020)
 \\   \hline
Case 3&mean & 0.466 &-1.06& 0.244 & 0.244 & 0.247 & 0.981
 \\
&s.d. & (0.005)&(1.152) & (0.012) & (0.012) & (0.012) & (0.047)
 \\   \hline
\end{tabular}
\end{center}
\end{table}

Table \ref{tab1} is the simulation results of 
the means and the standard deviations (s.d.s) of $\hat{\alpha}$, 
$\widehat{\theta}_0$, $\widehat{\theta}_1$, $\widehat{\eta}_1$, $\widehat{\theta}_2$ 
and $\widehat{\sigma}^2$.
In Case 1, $(m_1, m_2, b) = (200, 200, 0.005)$ for the estimation of $\alpha$ and $(n, m_1,m_2) = (100,30,30)$ for the estimation of $\theta_0$, $\theta_1$, 
$\eta_1$, $\theta_2$  and $\sigma^2$.
In Case 2, $(m_1, m_2, b) = (100, 100, 0.01)$ for the estimation of $\alpha$ and $(n, m_1,m_2) = (100,30,30)$ for the estimation of $\theta_0$, $\theta_1$, 
$\eta_1$, $\theta_2$  and $\sigma^2$.
In Case 3, $(m_1, m_2, b) = (50, 50, 0.02)$ for the estimation of $\alpha$ and $(n, m_1,m_2) = (100,30,30)$ for the estimation of $\theta_0$, $\theta_1$, 
$\eta_1$, $\theta_2$  and $\sigma^2$.

It seems from Table \ref{tab1} that the biases of $\widehat \alpha$ 
are small and get smaller as $m_1$ and $m_2$ increase. 
Moreover, Table \ref{tab1} indicates that in Case 1, the biases of 
$\widehat{\theta}_1$, $\widehat{\eta}_1$, $\widehat{\theta}_2$ and $\widehat{\sigma}^2$
are very small and $\widehat{\theta}_0$ has a small bias.

\section{Proofs}\label{sec6}
We set the following notation.
\begin{enumerate}
\item[1.]
For $a,b \in \mathbb R$, 
we write $a \lesssim b$ if $|a| \le C|b|$ for some constant $C>0$.

\item[2.]
For two functions $f, g$ such that $f, g : \mathbb R^d \to \mathbb R$,
we write $f(x) \lesssim g(x)$ $(x \to a)$ 
if $f(x) \lesssim g(x)$ in a neighborhood of $x = a$. 

\item[3.]
For $x=(x_1,\ldots, x_d) \in \mathbb R^d$ and a function $f$ such that $f:\mathbb R^d \to \mathbb R$,
we write $\pd_{x_i} f(x) = \frac{\pd}{\pd x_i}f(x)$,
$\pd_x f(x) = (\pd_{x_1}f(x), \ldots, \pd_{x_d}f(x))$ and 
$\pd_x^2 f(x) = (\pd_{x_j}\pd_{x_i}f(x))_{i,j=1}^d$.
\end{enumerate}

\subsection{Proof of Theorem \ref{th3-1}}
We provide only the proof of the $Q_1$-model.
For $\vartheta = (\kappa, \eta, \theta_2, \sigma^2)$, we define
\begin{equation}\label{pf1-1}
g_{r,\alpha}(\vartheta) = 
\frac{\sigma^2 \psi_{r,\alpha}(\theta_2)}{(1-2b)^2} 
\int_{b}^{1-b} \int_{b}^{1-b} \ee^{-\kappa y -\eta z} \dd y \dd z
\end{equation}
with $\psi_{r,\alpha}$ from \eqref{psi}. Since 
\begin{equation*}
J_0(\sqrt{2}x)-2J_0(x)+1
= \frac{2}{\pi} \int_0^{\pi/2} (1-\cos(x \cos(t)))(1-\cos(\sin(t))) \dd t \ge 0,
\quad x \ge 0,
\end{equation*}
we note that $g_{r,\alpha}(\vartheta) > 0$.

\textbf{Step 1:}
We first show
\begin{equation}\label{pf1-2}
\frac{1}{m N \Delta^{\alpha^*}} \sum_{k=1}^{m_2} \sum_{j=1}^{m_1} \sum_{i=1}^N 
\EE[(T_{i,j,k}X)^2]
= g_{r,\alpha^*}(\vartheta)+ \OO (\Delta).
\end{equation}
We see from Proposition 2.1 in \cite{TKU2024b} that
\begin{equation*}
\frac{1}{N \Delta^{\alpha^*}}\sum_{i=1}^N \EE[(T_{i,j,k}X)^2]
= \sigma^2 \ee^{-\kappa \overline y_j -\eta \overline z_k} \psi_{r,\alpha^*}(\theta_2)
+ \OO (\Delta).
\end{equation*}
Since it follows that for $f \in C^2([b,1-b])$,
\begin{align*}
&\Biggl| 
\frac{1}{m_1} \sum_{j=1}^{m_1} f(\overline y_j)
- \frac{1}{1-2b} \int_b^{1-b} f(y) \dd y
\Biggr|
\\
&=
\Biggl| \frac{1}{1-2b} 
\sum_{j=1}^{m_1} \int_{\widetilde y_{j-1}}^{\widetilde y_j} 
(f(\overline y_j) -f(y)) \dd y
\Biggr|
\\
&\le 
\Biggl| \frac{1}{1-2b} 
\sum_{j=1}^{m_1} \int_{\widetilde y_{j-1}}^{\widetilde y_j} 
f'(\overline y_j)(y -\overline y_j) \dd y
\Biggr|
\\
&\qquad +
\Biggl| \frac{1}{1-2b} 
\sum_{j=1}^{m_1} \int_{\widetilde y_{j-1}}^{\widetilde y_j} 
\biggl(
\int_0^1 f''(\overline y_j +u(y -\overline y_j)) \dd u 
\biggr)
(y -\overline y_j)^2 \dd y
\Biggr|
\\
&=
0 + \OO \Biggl( \sum_{j=1}^{m_1} 
\int_{\widetilde y_{j-1}}^{\widetilde y_j} (y -\overline y_j)^2 \dd y \Biggr)
\\
&= \OO(\delta^2) = \OO(\Delta),
\end{align*}
we obtain
\begin{align*}
\frac{1}{m N \Delta^{\alpha^*}} \sum_{k=1}^{m_2} \sum_{j=1}^{m_1} \sum_{i=1}^N 
\EE[(T_{i,j,k}X)^2]
&= \sigma^2 \psi_{r,\alpha^*}(\theta_2) \times
\frac{1}{m} \sum_{k=1}^{m_2} \sum_{j=1}^{m_1} 
\ee^{-\kappa \overline y_j -\eta \overline z_k}
+ \OO (\Delta)
\\
&= g_{r,\alpha^*}(\vartheta)+ \OO (\Delta).
\end{align*}

\textbf{Step 2:}
We next show 
\begin{equation}\label{pf1-3}
\sqrt{m N} \Biggl(
\frac{1}{m N \Delta^{\alpha^*}} \sum_{k=1}^{m_2} \sum_{j=1}^{m_1} \sum_{i=1}^N 
(T_{i,j,k}X)^2
-g_{r,\alpha^*}(\vartheta)
\Biggr) = \OO_\PP(1).
\end{equation}
We find from Lemma 4.11 in \cite{TKU2024b} that 
\begin{align*}
&\EE \Biggl[ \biggl(
\sum_{k=1}^{m_2} \sum_{j=1}^{m_1} \sum_{i=1}^N 
\Bigl( (T_{i,j,k}X)^2 -\EE[(T_{i,j,k}X)^2] \Bigr)
\biggr)^2 \Biggr]
\\
&= 
\sum_{k,k'=1}^{m_2} \sum_{j,j'=1}^{m_1} \sum_{i,i'=1}^{N}
\cov \Bigl[(T_{i,j,k}X)^2, (T_{i',j',k'}X)^2 \Bigr]
\\
&= \OO(m N \Delta^{2 \alpha^*}),
\end{align*}
which together with \eqref{pf1-2} yields
\begin{align*}
&\sqrt{m N} \Biggl(
\frac{1}{m N \Delta^{\alpha^*}} \sum_{k=1}^{m_2} \sum_{j=1}^{m_1} \sum_{i=1}^N 
(T_{i,j,k}X)^2 -g_{r,\alpha^*}(\vartheta)
\Biggr)
\\
&=
\frac{1}{\sqrt{m N} \Delta^{\alpha^*}} \sum_{k=1}^{m_2} \sum_{j=1}^{m_1} \sum_{i=1}^N 
\Bigl( (T_{i,j,k}X)^2 -\EE[(T_{i,j,k}X)^2] \Bigr)
\\
&\qquad +
\sqrt{m N} \Biggl(
\frac{1}{m N \Delta^{\alpha^*}} \sum_{k=1}^{m_2} \sum_{j=1}^{m_1} \sum_{i=1}^N 
\EE[(T_{i,j,k}X)^2] -g_{r,\alpha^*}(\vartheta)
\Biggr)
\\
&=
\OO_\PP(1).
\end{align*}

\textbf{Step 3:}
Finally, we prove 
$\sqrt{m N}(\widehat \alpha -\alpha^*) = \OO_\PP(1)$.
We define
\begin{equation*}
\mathcal Z =
\frac{1}{m N \Delta^{\alpha^*}} \sum_{k=1}^{m_2} \sum_{j=1}^{m_1} \sum_{i=1}^N 
(T_{i,j,k}X)^2,
\quad
\mathcal Z' =
\frac{1}{m' N' (\Delta')^{\alpha^*}} \sum_{k=1}^{m'_2} \sum_{j=1}^{m'_1} 
\sum_{i=1}^{N'} (T_{i,j,k}' X)^2.
\end{equation*}
In the same way as in Step 2, we have
\begin{equation}\label{pf1-4}
\sqrt{m N} \bigl( \mathcal Z' -g_{r,\alpha^*}(\vartheta) \bigr) = \OO_\PP(1).
\end{equation}
Since we have
\begin{equation*}
\biggl| \frac{\mathcal Z'}{\mathcal Z}-1 \biggr|
= \biggl| \frac{\mathcal Z'-g_{r,\alpha^*}(\vartheta) 
-(\mathcal Z -g_{r,\alpha^*}(\vartheta))}
{\mathcal Z -g_{r,\alpha^*}(\vartheta) +g_{r,\alpha^*}(\vartheta)} \biggr|
\le 
\frac{|\mathcal Z'- g_{r,\alpha^*}(\vartheta)| +|\mathcal Z -g_{r,\alpha^*}(\vartheta)|}
{|g_{r,\alpha^*}(\vartheta) -|\mathcal Z -g_{r,\alpha^*}(\vartheta)||},
\end{equation*}
we obtain by \eqref{pf1-3} and \eqref{pf1-4},
\begin{equation*}
\sqrt{m N} \biggl| \frac{\mathcal Z'}{\mathcal Z}-1 \biggr|
\le 
\frac{\sqrt{m N}|\mathcal Z'- g_{r,\alpha^*}(\vartheta)| 
+\sqrt{m N}|\mathcal Z -g_{r,\alpha^*}(\vartheta)|}
{|g_{r,\alpha^*}(\vartheta) -|\mathcal Z -g_{r,\alpha^*}(\vartheta)||}
= \OO_\PP(1).
\end{equation*}
Using the Taylor expansion
\begin{equation*}
\log(x) = (x-1) - \int_0^1 \frac{1-u}{(1+u(x-1))^2} \dd u (x-1)^2
\end{equation*}
and $\sqrt{m N} (\mathcal Z'/\mathcal Z -1) = \OO_\PP(1)$, we have
\begin{align*}
\sqrt{m N}(\widehat \alpha -\alpha^*)
&= \sqrt{m N} \left(
\log \left( \frac{\displaystyle\frac{1}{m' N'} \sum_{k=1}^{m'_2} \sum_{j=1}^{m'_1} 
\sum_{i=1}^{N'} (T_{i,j,k}'X)^2}
{\displaystyle\frac{1}{m N} \sum_{k=1}^{m_2} \sum_{j=1}^{m_1} \sum_{i=1}^N 
(T_{i,j,k}X)^2}
\right)/\log(4)
-\alpha^*
\right)
\\
&=
\sqrt{m N} \Biggl(
\log \biggl( \frac{(\Delta')^{\alpha^*} \mathcal Z'}
{\Delta^{\alpha^*} \mathcal Z} \biggr)
/\log(4)
-\alpha^*
\Biggr)
\\
&=
\sqrt{m N} \Biggl(
\biggl(
\log \biggl( \frac{\mathcal Z'}{\mathcal Z} \biggr)
+\alpha^* \log(4)
\biggr)/\log(4)
-\alpha^*
\Biggr)
\\
&=
\frac{\sqrt{m N}}{\log(4)}
\log \biggl( \frac{\mathcal Z'}{\mathcal Z} \biggr)
\\
&=
\frac{\sqrt{m N}}{\log(4)}
\Biggl( \frac{\mathcal Z'}{\mathcal Z} -1
+ \OO_\PP \biggl( \frac{1}{m N} \biggr)
\Biggr)
\\
&=\OO_\PP(1).
\end{align*}
This concludes the proof.

\subsection{Proofs of the results in Subsection \ref{sec4-1}}
\subsubsection{Proof of Theorem \ref{th4-1}}
Using the mean value theorem, we have
\begin{equation*}
-\frac{\sqrt{m N}}{\log(N)} 
\pd_{\vartheta} \mathcal K_{\mathbf{m},N} (\vartheta^*;\widehat \alpha)^\TT
=\int_0^1 \pd_{\vartheta}^2 
\mathcal K_{\mathbf{m},N} (\vartheta^* +u(\widehat\vartheta -\vartheta^*);\widehat \alpha) \dd u 
\frac{\sqrt{m N}}{\log(N)} (\widehat\vartheta-\vartheta^*).
\end{equation*}
We have already obtained the following results by \cite{TKU2024b}.
\begin{itemize}
\item[(i)]
There exists a function $\Xi \ni \vartheta \mapsto \mathcal K(\vartheta,\vartheta^*)$
which takes a unique minimum value at $\vartheta = \vartheta^*$ and satisfies
\begin{equation*}
\sup_{\vartheta \in \Xi} |\mathcal K_{\mathbf{m},N}(\vartheta;\alpha^*)
-\mathcal K(\vartheta,\vartheta^*)| 
= \oo_\PP(1),
\end{equation*}

\item[(ii)]
$\sqrt{m N} \pd_\vartheta \mathcal K_{\mathbf{m},N}(\vartheta^*;\alpha^*) = \OO_\PP(1)$,

\item[(iii)]
there exists a positive definite matrix $\mathcal V(\vartheta^*;\alpha^*)$ 
such that for $\epsilon_{\mathbf{m},N} \to 0$,
\begin{equation*}
\sup_{|\vartheta-\vartheta^*| \le \epsilon_{\mathbf{m},N}}
|\pd_{\vartheta}^2 \mathcal K_{\mathbf{m},N}(\vartheta;\alpha^*)
-\mathcal V(\vartheta^*;\alpha^*)|
= \oo_\PP(1).
\end{equation*}
\end{itemize}
Therefore, we will show the following three statements 
in order to obtain the desired results.
\begin{equation*}
\sup_{\vartheta \in \Xi} 
\bigl| \mathcal K_{\mathbf{m},N}(\vartheta;\widehat \alpha)
- \mathcal K_{\mathbf{m},N}(\vartheta;\alpha^*) \bigr| = \oo_\PP(1),
\end{equation*}
\begin{equation*}
\frac{\sqrt{m N}}{\log(N)} \bigl(
\pd_{\vartheta} \mathcal K_{\mathbf{m},N}(\vartheta^*;\widehat \alpha)
- \pd_{\vartheta} \mathcal K_{\mathbf{m},N}(\vartheta^*;\alpha^*) \bigr) = \OO_\PP(1),
\end{equation*}
\begin{equation*}
\sup_{\vartheta \in \Xi}
\bigl| \pd_{\vartheta}^2 \mathcal K_{\mathbf{m},N}(\vartheta;\widehat \alpha)
- \pd_{\vartheta}^2 \mathcal K_{\mathbf{m},N}(\vartheta;\alpha^*) \bigr| = \oo_\PP(1).
\end{equation*}

We write
$\mathcal K_{\mathbf{m},N}(\vartheta;\alpha) 
= \mathcal K_{\mathbf{m},N}^{(1)}(\vartheta;\alpha) +\mathcal K_{\mathbf{m},N}^{(2)}(\vartheta;\alpha)$
with 
\begin{align*}
\mathcal K_{\mathbf{m},N}^{(1)}(\vartheta;\alpha)
&= \frac{1}{m}\sum_{k=1}^{m_2}\sum_{j=1}^{m_1} 
\Biggl\{
\frac{1}{N \Delta^{\alpha}}\sum_{i=1}^{N} (T_{i,j,k}X)^2
-f_{r,\alpha}(\overline y_j, \overline z_k;\vartheta)
\Biggr\}^2,
\\
\mathcal K_{\mathbf{m},N}^{(2)}(\vartheta;\alpha)
&= \frac{1}{m} \sum_{k=1}^{m_2} \sum_{j=1}^{m_1} 
\Biggl\{
\frac{1}{N (2\Delta)^{\alpha}}\sum_{i=1}^{N-1} (\widetilde T_{i,j,k}X)^2
-f_{r/\sqrt{2},\alpha}(\overline y_j, \overline z_k;\vartheta)
\Biggr\}^2.
\end{align*}
We define
\begin{align*}
\mathbf{K}_{j,k}^{(1)}(\vartheta;\alpha)
&= \frac{1}{N \Delta^{\alpha}}\sum_{i=1}^{N} (T_{i,j,k}X)^2
-f_{r,\alpha}(\overline y_j, \overline z_k; \vartheta),
\\
\mathbf{K}_{j,k}^{(2)}(\vartheta;\alpha)
&= \frac{1}{N (2\Delta)^{\alpha}}\sum_{i=1}^{N-1} (\widetilde T_{i,j,k}X)^2
-f_{r/\sqrt{2},\alpha}(\overline y_j, \overline z_k; \vartheta),
\end{align*}
and $g_{j,k}^{l,p}(\vartheta;\alpha) 
= \pd_\vartheta^l f_{r/\sqrt{p},\alpha} (\overline y_j, \overline z_k; \vartheta)$
for $l= 0,1,2$ and $p = 1,2$.
Since
\begin{align*}
\mathcal K_{\mathbf{m},N}^{(p)}(\vartheta;\alpha)
&= \frac{1}{m} \sum_{k=1}^{m_2}\sum_{j=1}^{m_1} \mathbf{K}_{j,k}^{(p)}(\vartheta;\alpha)^2,
\\
\pd_{\vartheta} \mathcal K_{\mathbf{m},N}^{(p)}(\vartheta;\alpha)
&= -\frac{2}{m} \sum_{k=1}^{m_2}\sum_{j=1}^{m_1} 
\mathbf{K}_{j,k}^{(p)}(\vartheta;\alpha) g_{j,k}^{1,p}(\vartheta;\alpha),
\\
\pd_{\vartheta}^2 \mathcal K_{\mathbf{m},N}^{(p)}(\vartheta;\alpha)
&= \frac{2}{m} \sum_{k=1}^{m_2}\sum_{j=1}^{m_1} 
g_{j,k}^{1,p}(\vartheta;\alpha)
g_{j,k}^{1,p}(\vartheta;\alpha)^\TT
- \frac{2}{m} \sum_{k=1}^{m_2}\sum_{j=1}^{m_1}
\mathbf{K}_{j,k}^{(p)}(\vartheta;\alpha) g_{j,k}^{2,p}(\vartheta;\alpha)
\end{align*}
for $p = 1,2$, we have
\begin{align*}
&\Bigl| \mathcal K_{\mathbf{m},N}^{(p)}(\vartheta;\widehat \alpha) 
-\mathcal K_{\mathbf{m},N}^{(p)}(\vartheta;\alpha^*) \Bigr|
\\
& = \Biggl| \frac{1}{m}\sum_{k=1}^{m_2}\sum_{j=1}^{m_1} 
\mathbf{K}_{j,k}^{(p)}(\vartheta;\widehat \alpha)^2
- \frac{1}{m}\sum_{k=1}^{m_2}\sum_{j=1}^{m_1} 
\mathbf{K}_{j,k}^{(p)}(\vartheta;\alpha^*)^2
\Biggr|
\\
&=
\Biggl| \frac{1}{m}\sum_{k=1}^{m_2}\sum_{j=1}^{m_1} 
\bigl(
\mathbf{K}_{j,k}^{(p)}(\vartheta;\widehat \alpha) -\mathbf{K}_{j,k}^{(p)}(\vartheta;\alpha^*)
\bigr)^2
\\
&\qquad+
\frac{2}{m}\sum_{k=1}^{m_2}\sum_{j=1}^{m_1} 
\mathbf{K}_{j,k}^{(p)}(\vartheta;\alpha^*)
\bigl(
\mathbf{K}_{j,k}^{(p)}(\vartheta;\widehat \alpha) -\mathbf{K}_{j,k}^{(p)}(\vartheta;\alpha^*)
\bigr)
\Biggr|
\\
&\le
\frac{1}{m}\sum_{k=1}^{m_2}\sum_{j=1}^{m_1} 
\bigl(
\mathbf{K}_{j,k}^{(p)}(\vartheta;\widehat \alpha) -\mathbf{K}_{j,k}^{(p)}(\vartheta;\alpha^*)
\bigr)^2
\\
&\qquad+
2\Biggl\{\frac{1}{m}\sum_{k=1}^{m_2}\sum_{j=1}^{m_1} 
\mathbf{K}_{j,k}^{(p)}(\vartheta;\alpha^*)^2 \Biggr\}^{1/2}
\Biggl\{
\frac{1}{m}\sum_{k=1}^{m_2}\sum_{j=1}^{m_1} 
\bigl(
\mathbf{K}_{j,k}^{(p)}(\vartheta;\widehat \alpha) -\mathbf{K}_{j,k}^{(p)}(\vartheta;\alpha^*)
\bigr)^2 \Biggr\}^{1/2},
\end{align*}
\begin{align*}
&\Bigl| \pd_{\vartheta} \mathcal K_{\mathbf{m},N}^{(1)}(\vartheta;\widehat \alpha)
-\pd_{\vartheta} \mathcal K_{\mathbf{m},N}^{(1)}(\vartheta;\alpha^*) \Bigr|
\\
&= 2 \Biggl|
\frac{1}{m} \sum_{k=1}^{m_2} \sum_{j=1}^{m_1}
\Bigl(
\bigl(\mathbf{K}_{j,k}^{(p)}(\vartheta;\widehat \alpha)
-\mathbf{K}_{j,k}^{(p)}(\vartheta;\alpha^*) \bigr)
g_{j,k}^{1,p}(\vartheta;\widehat \alpha)
\\
&\qquad \qquad +\mathbf{K}_{j,k}^{(p)}(\vartheta;\alpha^*)
\bigl( g_{j,k}^{1,p}(\vartheta;\widehat \alpha)- g_{j,k}^{1,p}(\vartheta;\alpha^*) \bigr) 
\Bigr)
\Biggr|
\\
&\le 2 
\Biggl\{
\frac{1}{m} \sum_{k=1}^{m_2} \sum_{j=1}^{m_1}
\bigl( \mathbf{K}_{j,k}^{(p)}(\vartheta;\widehat \alpha)
-\mathbf{K}_{j,k}^{(p)}(\vartheta;\alpha^*) \bigr)^2
\Biggr\}^{1/2}
\Biggl\{
\frac{1}{m} \sum_{k=1}^{m_2} \sum_{j=1}^{m_1}
|g_{j,k}^{1,p}(\vartheta;\widehat \alpha)|^2
\Biggr\}^{1/2}
\\
&\qquad +
2 \Biggl\{
\frac{1}{m} \sum_{k=1}^{m_2} \sum_{j=1}^{m_1}
\mathbf{K}_{j,k}^{(p)}(\vartheta;\alpha^*)^2
\Biggr\}^{1/2}
\Biggl\{
\frac{1}{m} \sum_{k=1}^{m_2} \sum_{j=1}^{m_1}
\bigl| g_{j,k}^{1,p}(\vartheta;\widehat \alpha)
-g_{j,k}^{1,p}(\vartheta;\alpha^*) \bigr|^2
\Biggr\}^{1/2},
\end{align*}
\begin{align*}
&\Bigl| \pd_{\vartheta}^2 \mathcal K_{\mathbf{m},N}^{(1)}(\vartheta;\widehat \alpha)
-\pd_{\vartheta}^2 \mathcal K_{\mathbf{m},N}^{(1)}(\vartheta;\alpha^*) \Bigr|
\\
&\le 
\frac{2}{m} \sum_{k=1}^{m_2} \sum_{j=1}^{m_1}
\bigl|
g_{j,k}^{1,p}(\vartheta;\widehat \alpha) g_{j,k}^{1,p}(\vartheta;\widehat \alpha)^\TT
-g_{j,k}^{1,p}(\vartheta;\alpha^*) g_{j,k}^{1,p}(\vartheta;\alpha^*)^\TT
\bigr|
\\
&\qquad +2\Biggl\{
\frac{1}{m} \sum_{k=1}^{m_2} \sum_{j=1}^{m_1}
\bigl( \mathbf{K}_{j,k}^{(p)}(\vartheta;\widehat \alpha)
-\mathbf{K}_{j,k}^{(p)}(\vartheta;\alpha^*) \bigr)^2
\Biggr\}^{1/2}
\Biggl\{
\frac{1}{m} \sum_{k=1}^{m_2} \sum_{j=1}^{m_1}
|g_{j,k}^{2,p}(\vartheta;\widehat \alpha)|^2
\Biggr\}^{1/2}
\\
&\qquad +
2 \Biggl\{
\frac{1}{m} \sum_{k=1}^{m_2} \sum_{j=1}^{m_1}
\mathbf{K}_{j,k}^{(p)}(\vartheta;\alpha^*)^2
\Biggr\}^{1/2}
\Biggl\{
\frac{1}{m} \sum_{k=1}^{m_2} \sum_{j=1}^{m_1}
\bigl| g_{j,k}^{2,p}(\vartheta;\widehat \alpha) 
-g_{j,k}^{2,p}(\vartheta;\alpha^*) \bigr|^2
\Biggr\}^{1/2}.
\end{align*}
Therefore, we will show that for $p = 1,2$,
\begin{equation}\label{pf4-1-1}
\frac{1}{m} \sup_{\vartheta \in \Xi} \sum_{k=1}^{m_2} \sum_{j=1}^{m_1}
\bigl| g_{j,k}^{l,p}(\vartheta;\widehat \alpha) 
-g_{j,k}^{l,p}(\vartheta;\alpha^*) \bigr|^2
= \OO_\PP \biggl(\frac{1}{m N} \biggr),
\quad l = 0, 1, 2,
\end{equation}
\begin{equation}\label{pf4-1-2}
\frac{1}{m} \sup_{\vartheta \in \Xi}
\sum_{k=1}^{m_2}\sum_{j=1}^{m_1} 
\mathbf{K}_{j,k}^{(p)}(\vartheta;\alpha^*)^2
=\OO_\PP(1),
\end{equation}
\begin{equation}\label{pf4-1-3}
\frac{1}{m} \sup_{\vartheta \in \Xi}
\sum_{k=1}^{m_2}\sum_{j=1}^{m_1} 
\bigl(
\mathbf{K}_{j,k}^{(p)}(\vartheta;\widehat \alpha) -\mathbf{K}_{j,k}^{(p)}(\vartheta;\alpha^*)
\bigr)^2
= \OO_\PP \bigl( (\Delta \log (N))^2 \bigr).
\end{equation}

\textit{Proof of \eqref{pf4-1-1}. }
For any $\theta_2 > 0$, $r >0$ and $l=0,1,2$, 
$\pd_{\theta_2}^l \psi_{r,\alpha}(\theta_2)$ 
is differentiable with respect to $\alpha \in (\underline \alpha, \overline \alpha)$.
Indeed, we express
\begin{equation*}
\psi_{r,\alpha} (\theta_2)
= \frac{2}{\theta_2^{1+\alpha} \pi} 
\int_0^\infty \frac{1-\ee^{-\theta_2 x^2}}{x^{1+2\alpha}} 
\bigl( J_0(\sqrt{2}r x) -2J_0(r x)+1 \bigr) \dd x
\end{equation*}
and have
\begin{align*}
&\int_0^\infty 
\frac{\pd}{\pd \alpha}
\biggl(
\frac{1-\ee^{-\theta_2 x^2}}{x^{1+2\alpha}} 
\bigl( J_0(\sqrt{2}r x) -2J_0(r x)+1 \bigr) \biggr) \dd x
\\
&= 2 \int_0^\infty 
\frac{1-\ee^{-\theta_2 x^2}}{x^{1+2\alpha}} \log \biggl( \frac{1}{x} \biggr)
\bigl( J_0(\sqrt{2}r x) -2J_0(r x)+1 \bigr) \dd x.
\end{align*}
Since it holds that for any $\alpha \in (\underline \alpha, \overline \alpha)$,
\begin{equation*}
(0<) \int_0^\infty \frac{1-\ee^{-\theta_2 x^2}}{x^{1+2\alpha}} 
\bigl( J_0(\sqrt{2}r x) -2J_0(r x)+1 \bigr) \dd x < \infty
\end{equation*}
and for any $a >0$, 
\begin{equation*}
\lim_{x \searrow 0} x^a \log(x) = 0,
\quad
\lim_{x \to \infty} \frac{\log(x)}{x^a} = 0,
\end{equation*}
we obtain
\begin{equation*}
\int_0^\infty 
\biggl|
\frac{1-\ee^{-\theta_2 x^2}}{x^{1+2\alpha}} \log \biggl( \frac{1}{x} \biggr)
\bigl( J_0(\sqrt{2}r x) -2J_0(r x)+1 \bigr) \biggr|\dd x
< \infty
\end{equation*}
and see $\psi_{r,\alpha}(\theta_2)$ is differentiable with respect to $\alpha$. 
It can be shown for $l = 1, 2$ in the same way.
There exist a closed interval 
$I_{\alpha^*} \subset (\underline \alpha, \overline \alpha)$ 
and $\epsilon > 0$ such that 
$\alpha_u = \alpha^* +u(\widehat \alpha -\alpha^*) \in I_{\alpha^*}$
for any $u \in [0,1]$ on 
$\Omega_\epsilon =  \{ |\widehat \alpha -\alpha^*| < \epsilon \}$. 
It then follows that on $\Omega_\epsilon$,
\begin{align*}
\bigl| \pd_{\theta_2}^l \psi_{r,\widehat \alpha}(\theta_2) 
-\pd_{\theta_2}^l \psi_{r,\alpha^*}(\theta_2) \bigr|
&= \biggl| 
\int_0^1 \pd_\alpha \pd_{\theta_2}^l 
\psi_{r,\alpha_u}(\theta_2) \dd u
\biggr|
|\widehat \alpha -\alpha^*|
\\
&\le \frac{1}{\sqrt{m N}} \times 
\sup_{\alpha \in I_{\alpha^*}}
\bigl| \pd_\alpha \pd_{\theta_2}^l \psi_{r,\alpha}(\theta_2) \bigr|
\sqrt{m N}|\widehat \alpha -\alpha^*|
\end{align*}
for sufficient large $m$, $N$ and $l=0,1,2$. Since
\begin{equation*}
f_{r,\widehat \alpha}(y, z, \vartheta)
-f_{r,\alpha^*}(y, z, \vartheta)
= \sigma^2 \ee^{-\kappa y -\eta z} 
(\psi_{r,\widehat \alpha}(\theta_2) -\psi_{r,\alpha^*}(\theta_2)),
\end{equation*}
\begin{equation*}
\sup_{\vartheta \in \Xi} \bigl| \pd_\vartheta^l f_{r,\widehat \alpha}(y, z, \vartheta)
-\pd_\vartheta^l f_{r,\alpha^*}(y, z, \vartheta) \bigr|
\lesssim
\sum_{\nu=0}^l \sup_{\theta_2} 
\bigl| \pd_{\theta_2}^\nu \psi_{r,\widehat \alpha}(\theta_2) 
-\pd_{\theta_2}^\nu \psi_{r,\alpha^*}(\theta_2) \bigr|, 
\quad l=0,1,2
\end{equation*}
the continuity of the function
$(\theta_2, \alpha) \mapsto \pd_\alpha \pd_{\theta_2}^l \psi_{r,\alpha}(\theta_2)$ 
for $l=0,1,2$ and Theorem \ref{th3-1}, we get the desired result.
 
\textit{Proof of \eqref{pf4-1-2}. }
Since $\widetilde T_{i,j,k} X = T_{i,j,k}X +T_{i+1,j,k}X$ 
and $2^{-\alpha} \le 1$ for $\alpha > 0$, we have 
\begin{equation}\label{eq1-pf4-1-2}
\frac{1}{N (2\Delta)^{\alpha}}\sum_{i=1}^{N-1} (\widetilde T_{i,j,k}X)^2
\lesssim 
\frac{1}{N \Delta^{\alpha}}\sum_{i=1}^{N} (T_{i,j,k}X)^2
\end{equation}
and 
\begin{equation*}
\frac{1}{m}\sum_{k=1}^{m_2}\sum_{j=1}^{m_1} 
\mathbf{K}_{j,k}^{(p)}(\vartheta;\alpha)^2
\lesssim
\frac{1}{m}\sum_{k=1}^{m_2}\sum_{j=1}^{m_1}
\biggl(
\frac{1}{N \Delta^{\alpha}}\sum_{i=1}^{N} (T_{i,j,k}X)^2
\biggr)^2
+ \frac{1}{m}\sum_{k=1}^{m_2}\sum_{j=1}^{m_1} g_{j,k}^{0,p}(\vartheta;\alpha)^2.
\end{equation*}
By Proposition 2.1 and Lemma 4.11 in \cite{TKU2024b}, we obtain
\begin{align*}
&\EE \Biggl[
\biggl( \frac{1}{N \Delta^{\alpha^*}}\sum_{i=1}^{N} (T_{i,j,k}X)^2 \biggr)^2
\Biggr]
\\
&\lesssim \EE \Biggl[
\biggl( \frac{1}{N \Delta^{\alpha^*}}\sum_{i=1}^{N} \bigl( (T_{i,j,k}X)^2 
-\EE[(T_{i,j,k}X)^2] \bigr)
\biggr)^2
\Biggr]
+ \biggl( \frac{1}{N \Delta^{\alpha^*}}\sum_{i=1}^{N} 
\EE[(T_{i,j,k}X)^2]
\biggr)^2
\\
&= \VV \Biggl[\frac{1}{N \Delta^{\alpha^*}}\sum_{i=1}^{N} (T_{i,j,k}X)^2 \Biggr]
+\OO(1)
\\
&= \frac{1}{(N \Delta^{\alpha^*})^2} \sum_{i,i'=1}^{N} 
\cov \bigl[(T_{i,j,k}X)^2, (T_{i',j,k}X)^2 \bigr]
+\OO(1)
\\
&= \OO(N^{-1}) +\OO(1)
=\OO(1),
\end{align*}
which together with 
\begin{equation*}
\sup_{\vartheta \in \Xi} 
\Biggl| \frac{1}{m}\sum_{k=1}^{m_2}\sum_{j=1}^{m_1} 
g_{j,k}^{0,p}(\vartheta;\alpha^*)^2
- \frac{1}{(1-2b)^2} \iint_{[b,1-b]^2} 
f_{r/\sqrt{p},\alpha^*}(y,z,\vartheta)^2 \dd y \dd z \Biggr|
\to 0
\end{equation*}
yields the desired result.

\textit{Proof of \eqref{pf4-1-3}. }
It follows from \eqref{eq1-pf4-1-2} and the Taylor expansion that 
\begin{align*}
&\Bigl| \mathbf{K}_{j,k}^{(p)}(\vartheta;\widehat \alpha) 
-\mathbf{K}_{j,k}^{(p)}(\vartheta;\alpha^*) \Bigr|
\\
&\lesssim
\Delta^{\alpha^*} 
\biggl| \frac{1}{\Delta^{\widehat \alpha}} -\frac{1}{\Delta^{\alpha^*}} \biggr|
\frac{1}{N \Delta^{\alpha^*}}\sum_{i=1}^{N} (T_{i,j,k}X)^2
+\bigl| g_{j,k}^{0,p}(\widehat \alpha; \vartheta)
-g_{j,k}^{0,p}(\alpha^*; \vartheta) \bigr|,
\end{align*}
and
\begin{align*}
\biggl| \frac{1}{\Delta^{\widehat \alpha}} -\frac{1}{\Delta^{\alpha^*}} \biggr|
&= \frac{1}{\Delta^{\alpha^*}} \log \biggl( \frac{1}{\Delta} \biggr)
\int_0^1 \biggl( \frac{1}{\Delta} \biggr)^{u(\widehat \alpha -\alpha^*)} \dd u
|\widehat \alpha -\alpha^*|
\\
&\le 
\frac{1}{\Delta^{\alpha^*}} \log \biggl( \frac{1}{\Delta} \biggr)
\biggl( \frac{1}{\Delta} \biggr)^{|\widehat \alpha -\alpha^*|}
|\widehat \alpha -\alpha^*|
\\
&= 
\frac{1}{\Delta^{\alpha^*}} \log \biggl( \frac{1}{\Delta} \biggr)
\exp \Biggl(|\widehat \alpha -\alpha^*| \log \biggl( \frac{1}{\Delta} \biggr) \Biggr)
|\widehat \alpha -\alpha^*|
\\
&= 
\frac{1}{\Delta^{\alpha^*}\sqrt{m N}} \log \biggl( \frac{1}{\Delta} \biggr)
\exp \Biggl(\sqrt{m N}|\widehat \alpha -\alpha^*| 
\times \frac{\log (N)}{\sqrt{m N}} \Biggr)
\sqrt{m N} |\widehat \alpha -\alpha^*|
\\
&= \OO_\PP \Biggl(\Delta^{1-\alpha^*} \log \biggl( \frac{1}{\Delta} \biggr) \Biggr).
\end{align*}
Therefore, we obtain
\begin{align*}
&\frac{1}{m} \sup_{\vartheta \in \Xi} \sum_{k=1}^{m_2}\sum_{j=1}^{m_1} 
\bigl(
\mathbf{K}_{j,k}^{(p)}(\vartheta;\widehat \alpha) -\mathbf{K}_{j,k}^{(p)}(\vartheta;\alpha^*)
\bigr)^2
\\
&\lesssim
\Delta^{2\alpha^*}
\biggl( \frac{1}{\Delta^{\widehat \alpha}} -\frac{1}{\Delta^{\alpha^*}} \biggr)^2
\times
\frac{1}{m}\sum_{k=1}^{m_2}\sum_{j=1}^{m_1} 
\biggl( \frac{1}{N \Delta^{\alpha^*}}\sum_{i=1}^{N} (T_{i,j,k}X)^2 \biggr)^2
\\
&\qquad+
\frac{1}{m} \sup_{\vartheta \in \Xi} \sum_{k=1}^{m_2}\sum_{j=1}^{m_1} 
\bigl( g_{j,k}^{0,p}(\widehat \alpha; \vartheta)
-g_{j,k}^{0,p}(\alpha^*; \vartheta) \bigr)^2
\\
&= \OO_\PP \bigl( (\Delta \log (N))^2 \bigr) +\OO_\PP \biggl(\frac{1}{m N} \biggr)
\\
&= \OO_\PP \bigl( (\Delta \log (N))^2 \bigr).
\end{align*}

\subsubsection{Proof of Theorem \ref{th4-2}}

\textit{Proof of \eqref{eq-th2}. }
First, we show that
\begin{equation}\label{eq-3-1}
\sqrt{n}
\begin{pmatrix}
\widetilde \sigma_{1,1}^2-(\sigma_{1,1}^*)^2
\\
\widetilde \sigma_{1,2}^2-(\sigma_{1,2}^*)^2
\end{pmatrix}
\dto 
N
\Biggl(0,2
\begin{pmatrix}
(\sigma_{1,1}^*)^4 & 0
\\
0 & (\sigma_{1,2}^*)^4 
\end{pmatrix}
\Biggr).
\end{equation}
In order to show \eqref{eq-3-1}, we verify that
under $\frac{n}{(M_1 \land M_2)^{\alpha^* \tand 1}} \to 0$, 
\begin{align*}
\sqrt{n}\Biggl(
\widetilde\sigma_{l_1,l_2}^2
-\sum_{i=1}^{n}
(x_{l_1,l_2}(\widetilde t_i)-x_{l_1,l_2}(\widetilde t_{i-1}))^2
\Biggr)
=\oo_\PP(1),
\end{align*}
that is, for $k = 1,2$,
\begin{equation*}
n \sum_{i=1}^n \mathcal B_{k,i}^2 = \oo_\PP(1),
\end{equation*}
where
$\widetilde \Delta_i X(y,z)=X_{\widetilde t_i}(y,z)-X_{\widetilde t_{i-1}}(y,z)$,
\begin{align*}
\mathcal B_{1,i} 
&= \sum_{j=1}^{M_1}\sum_{k=1}^{M_2}
\widetilde \Delta_i X(y_{j-1},z_{k-1})
(\delta_j^{[y]} g_{l_1}(\widehat \kappa) \delta_k^{[z]} g_{l_2}(\widehat \eta)
-\delta_j^{[y]} g_{l_1}(\kappa^*) \delta_k^{[z]} g_{l_2}(\eta^*)),
\\
\mathcal B_{2,i} &= 
\sum_{j=1}^{M_1}\sum_{k=1}^{M_2}
\int_{z_{k-1}}^{z_k} \int_{y_{j-1}}^{y_j}
(\widetilde \Delta_i X(y_{j-1},z_{k-1})-\widetilde \Delta_i X(y,z))
\sin(\pi l_1 y)\sin(\pi l_2 z)
\ee^{(\kappa^* y+\eta^* z)/2}
\dd y\dd z.
\end{align*}

Using the Taylor expansion, we have
\begin{equation*}
n \sum_{i=1}^n \mathcal B_{1,i}^2
\le\mathcal C_n
\times \frac{m N}{(\log(N))^2} 
(|\widehat \kappa -\kappa^*|^2 +|\widehat \eta -\eta^*|^2),
\end{equation*}
where $\mathcal C_n$ satisfies for $\epsilon_1 > 0$, 
\begin{align*}
\mathcal C_n
\lesssim
\frac{n (\log(N))^2}{mN}\sum_{i=1}^{n}
\frac{1}{M}
\sum_{j=1}^{M_1}\sum_{k=1}^{M_2}
(\widetilde \Delta_i X)^2(y_{j-1},z_{k-1})
\end{align*}
on $\Omega_1 = \{|\widehat \kappa -\kappa^*| +|\widehat \eta -\eta^*| < \epsilon_1 \}$.
Note that
\begin{align*}
\sum_{i=1}^n \EE \Bigl[ (\widetilde \Delta_i X)^2(y,z) \Bigr]
\lesssim n \sum_{l_1,l_2 \ge 1} \frac{1-\ee^{-\lambda_{l_1,l_2}\Delta_n}}
{\lambda_{l_1,l_2}^{1+\alpha^*}}
\end{align*}
uniformly in $y,z \in [0,1]$ and 
\begin{equation*}
\sum_{l_1,l_2 \ge 1} \frac{1-\ee^{-\lambda_{l_1,l_2}\Delta_n}}
{\lambda_{l_1,l_2}^{1+\alpha^*}}
= \OO(n^{-(\alpha^* \tand 1)}).
\end{equation*}
It then follows from 
$\frac{n^{2-(\alpha^* \tand 1)} (\log(N))^2}{m N} \to 0$ and Theorem \ref{th4-1} that
for $\epsilon_2 > 0$, 
\begin{equation*}
\PP(|\mathcal C_n|>\epsilon_2)
\lesssim
\frac{n^{2-(\alpha^* \tand 1)} (\log(N))^2}{\epsilon_2 m N} +\PP(\Omega_1^{\mathrm c})
\to 0.
\end{equation*}
Therefore, we have $n \sum_{i=1}^{n} \mathcal B_{1,i}^2 = \oo_\PP(1)$.

From Lemma 9 in \cite{TKU2024arXiv3}, we have
\begin{align*}
&\sup_{(y,z) \in (y_{j-1},y_j]\times (z_{k-1},z_k]}
\EE \Bigl[(\widetilde \Delta_i X(y_{j-1},z_{k-1})
-\widetilde \Delta_i X(y,z))^2 \Bigr]
\\
& \lesssim
\sum_{l_1, l_2 \ge 1}
\frac{1-\ee^{-\lambda_{l_1, l_2} \Delta_n}}{\lambda_{l_1, l_2}^{1+\alpha^*}}
\sup_{(y,z) \in (y_{j-1},y_j]\times (z_{k-1},z_k]}
(e_{l_1,l_2}(y_{j-1},z_{k-1}) -e_{l_1,l_2}(y,z))^2
\\
& = \OO \biggl(\frac{1}{(M_1 \land M_2)^{2(\alpha^* \tand 1)}} \biggr)
\end{align*}
uniformly in $i,j,k$.
Therefore, we see that
under $\frac{n^2}{(M_1 \land M_2)^{2(\alpha^* \tand 1)}} \to 0$, 
\begin{align*}
n\sum_{i=1}^{n}
\mathcal B_{2,i}^2
&\lesssim
n\sum_{i=1}^{n}
\sum_{j=1}^{M_1}\sum_{k=1}^{M_2}
\int_{z_{k-1}}^{z_k} \int_{y_{j-1}}^{y_j}
\bigl(\widetilde \Delta_i X(y_{j-1},z_{k-1})-\widetilde \Delta_i X(y,z)\bigr)^2
\dd y\dd z
\\
& = \OO_\PP \biggl( \frac{n^2}{(M_1 \land M_2)^{2 (\alpha^* \tand 1)}} \biggr)
= \oo_\PP(1).
\end{align*}

We next show \eqref{eq-th2}. 
It follows from Theorem \ref{th4-1} that for $g \in C^1(\Xi)$,
\begin{equation*}
\frac{\sqrt{m N}}{\log(N)}(g(\widehat \vartheta) -g(\vartheta^*)) = \OO_\PP(1)
\end{equation*}
and 
\begin{align*}
\sqrt{n} (\widetilde \theta_0 -\theta_0^*)
&= - \sqrt{n} (\widetilde \lambda_{1,1} -\lambda_{1,1}^* )
+\sqrt{n} \biggl\{
\biggl( \frac{\widehat \kappa^2+\widehat \eta^2}{4} +2\pi^2 \biggr) \widehat \theta_2
-\biggl( \frac{(\kappa^*)^2 +(\eta^*)^2}{4} +2\pi^2 \biggr) \theta_2^*
\biggr\}
\\
&= - \sqrt{n} (\widetilde \lambda_{1,1} -\lambda_{1,1}^* )
\\
&\qquad 
+\frac{\sqrt{n} \log(N)}{\sqrt{m N}} 
\frac{\sqrt{m N}}{\log(N)}
\biggl\{
\biggl( \frac{\widehat \kappa^2+\widehat \eta^2}{4} +2\pi^2 \biggr) \widehat \theta_2
-\biggl( \frac{(\kappa^*)^2 +(\eta^*)^2}{4} +2\pi^2 \biggr) \theta_2^*
\biggr\}
\\
&= - \sqrt{n} (\widetilde \lambda_{1,1} -\lambda_{1,1}^* ) +\oo_\PP(1).
\end{align*}
We decompose $\sqrt{n} (\widetilde \lambda_{1,1} -\lambda_{1,1}^* )$ as follows.
\begin{align*}
\sqrt{n} (\widetilde \lambda_{1,1} -\lambda_{1,1}^*)
&= \sqrt{n} 
\Biggl(
\biggl( \frac{\widehat \sigma^2}{\widetilde \sigma_{1,1}^2} \biggr)^{1/\widehat \alpha}
-\biggl( \frac{(\sigma^*)^2}{(\sigma_{1,1}^*)^2} \biggr)^{1/\alpha^*}
\Biggr)
\\
&= \sqrt{n} 
\Biggl(
\biggl( \frac{\widehat \sigma^2}{\widetilde \sigma_{1,1}^2} \biggr)^{1/\widehat \alpha}
-\biggl( \frac{\widehat \sigma^2}{\widetilde \sigma_{1,1}^2} \biggr)^{1/\alpha^*}
\Biggr)
\\
&\qquad +\sqrt{n} 
\Biggl(
\biggl( \frac{\widehat \sigma^2}{\widetilde \sigma_{1,1}^2} \biggr)^{1/\alpha^*}
-\biggl( \frac{(\sigma^*)^2}{\widetilde \sigma_{1,1}^2} \biggr)^{1/\alpha^*}
\Biggr)
\\
&\qquad +\sqrt{n} 
\Biggl(
\biggl( \frac{(\sigma^*)^2}{\widetilde \sigma_{1,1}^2} \biggr)^{1/\alpha^*}
-\biggl( \frac{(\sigma^*)^2}{(\sigma_{1,1}^*)^2} \biggr)^{1/\alpha^*}
\Biggr)
\\
&=: T_1 +T_2 +T_3.
\end{align*}
It then holds from the delta method that
\begin{align*}
T_1 
&= \sqrt{n} 
\Biggl(
\biggl( \frac{\widehat \sigma^2}{\widetilde \sigma_{1,1}^2} 
\biggr)^{1/\widehat \alpha -1/\alpha^*}
-1 \Biggr)
\biggl( \frac{\widehat \sigma^2}{\widetilde \sigma_{1,1}^2} \biggr)^{1/\alpha^*}
\\
&= \sqrt{n} 
\Biggl(
\exp \Biggl(
\biggl( \frac{1}{\widehat \alpha} -\frac{1}{\alpha^*} \biggr) 
\log \biggl( \frac{\widehat \sigma^2}{\widetilde \sigma_{1,1}^2} 
\biggr) \Biggr)
-1 \Biggr)
\biggl( \frac{\widehat \sigma^2}{\widetilde \sigma_{1,1}^2} \biggr)^{1/\alpha^*}
\\
&= \frac{\sqrt{n}}{\sqrt{m N}} 
\sqrt{m N}\biggl( \frac{1}{\widehat \alpha} -\frac{1}{\alpha^*} \biggr) 
\log \biggl( \frac{\widehat \sigma^2}{\widetilde \sigma_{1,1}^2} \biggr)
\\
&\qquad \times
\int_0^1 \exp \Biggl(
u \biggl( \frac{1}{\widehat \alpha} -\frac{1}{\alpha^*} \biggr) 
\log \biggl( \frac{\widehat \sigma^2}{\widetilde \sigma_{1,1}^2} \biggr) \Biggr) \dd u
\biggl( \frac{\widehat \sigma^2}{\widetilde \sigma_{1,1}^2} \biggr)^{1/\alpha^*}
\\
&= \oo_\PP(1).
\end{align*}
In the same way as the proof of Theorem 2.3 in \cite{TKU2024b}, we have
\begin{align*}
T_2 
&= \sqrt{n} (\widehat \sigma^{2/\alpha^*}-(\sigma^*)^{2/\alpha^*})
\widetilde \sigma_{1,1}^{-2/\alpha^*}
\\
&= 
\frac{\sqrt{n}\log(N)}{\sqrt{m N}}
\frac{\sqrt{m N}}{\log(N)} 
(\widehat \sigma^{2/\alpha^*}-(\sigma^*)^{2/\alpha^*})
\widetilde \sigma_{1,1}^{-2/\alpha^*}
\\
&=\oo_\PP(1),
\\
T_3 
&= (\sigma^*)^{2/\alpha^*} 
\sqrt{n} \Bigl( \widetilde \sigma_{1,1}^{-2/\alpha^*} 
-(\sigma_{1,1}^*)^{-2/\alpha^*} \Bigr)
\\
&= -\frac{(\sigma^*)^{2/\alpha^*}}{\alpha^*} 
(\sigma_{1,1}^*)^{-2/\alpha^*-2}
\sqrt{n} \bigl( \widetilde \sigma_{1,1}^2 -(\sigma_{1,1}^*)^2 \bigr)
+\oo_\PP(1)
\\
&= -\frac{(\sigma^*)^{-2}}{\alpha^*}
\bigl( (\lambda_{1,1}^*)^{1+\alpha^*}, 0 \bigr)
\sqrt{n} 
\begin{pmatrix}
\widetilde \sigma_{1,1}^2 -(\sigma_{1,1}^*)^2
\\
\widetilde \sigma_{1,2}^2 -(\sigma_{1,2}^*)^2
\end{pmatrix}
+\oo_\PP(1).
\end{align*}
Hence, we have
\begin{equation*}
\sqrt{n} (\widetilde \theta_0 -\theta_0^*)
= \frac{(\sigma^*)^{-2}}{\alpha^*}
\bigl( (\lambda_{1,1}^*)^{1+\alpha^*}, 0 \bigr)
\sqrt{n} 
\begin{pmatrix}
\widetilde \sigma_{1,1}^2 -(\sigma_{1,1}^*)^2
\\
\widetilde \sigma_{1,2}^2 -(\sigma_{1,2}^*)^2
\end{pmatrix}
+\oo_\PP(1).
\end{equation*}
Noting that
\begin{equation*}
\widetilde \theta_2 = 
\frac{\widetilde \lambda_{1,2}-\widetilde \lambda_{1,1}}{3\pi^2}
= \frac{\widehat \sigma^{2/\widehat \alpha}}{3\pi^2}
\bigl( \widetilde \sigma_{1,2}^{-2/\widehat \alpha} 
- \widetilde \sigma_{1,1}^{-2/\widehat \alpha} \bigr),
\end{equation*}
we see from the delta method that
\begin{align*}
\sqrt{n} (\widetilde \theta_2 - \theta_2^*) 
&= 
\frac{\sqrt{n}}{3\pi^2}
\Bigl(
\widehat \sigma^{2/\widehat \alpha} 
\bigl( \widetilde \sigma_{1,2}^{-2/\widehat \alpha} 
-\widetilde \sigma_{1,1}^{-2/\widehat \alpha} \bigr)
-(\sigma^*)^{2/\alpha^*} 
\bigl( (\sigma_{1,2}^*)^{-2/\alpha^*} - (\sigma_{1,1}^*)^{-2/\alpha^*} \bigr)
\Bigr)
\\
&=
\frac{\sqrt{n}}{3\pi^2}
\Bigl(
\widehat \sigma^{2/\widehat \alpha} 
\bigl( \widetilde \sigma_{1,2}^{-2/\widehat \alpha} 
-\widetilde \sigma_{1,1}^{-2/\widehat \alpha} \bigr)
-\widehat \sigma^{2/\alpha^*} 
\bigl( \widetilde \sigma_{1,2}^{-2/\alpha^*} 
-\widetilde \sigma_{1,1}^{-2/\alpha^*} \bigr)
\Bigr)
\\
&\qquad
+\frac{\sqrt{n}}{3\pi^2}
\bigl( \widehat \sigma^{2/\alpha^*} - (\sigma^*)^{2/\alpha^*} \bigr)
\bigl( \widetilde \sigma_{1,2}^{-2/\alpha^*} 
-\widetilde \sigma_{1,1}^{-2/\alpha^*} \bigr)
\\
&\qquad+ 
\frac{\sqrt{n}(\sigma^*)^{2/\alpha^*}}{3\pi^2}
\Bigl( 
\bigl( \widetilde \sigma_{1,2}^{-2/\alpha^*} 
-\widetilde \sigma_{1,1}^{-2/\alpha^*} \bigr)
- \bigl( (\sigma_{1,2}^*)^{-2/\alpha^*} - (\sigma_{1,1}^*)^{-2/\alpha^*} \bigr)
\Bigr)
\\
&=
\frac{\sqrt{n}}{3\pi^2}
\Biggl\{
\biggl( \frac{\widehat \sigma^2}{\widetilde \sigma_{1,2}^2} \biggr)^{1/\widehat \alpha} 
-\biggl( \frac{\widehat \sigma^2}{\widetilde \sigma_{1,2}^2} \biggr)^{1/\alpha^*}
-\Biggl(
\biggl( \frac{\widehat \sigma^2}{\widetilde \sigma_{1,1}^2} \biggr)^{1/\widehat \alpha} 
-\biggl( \frac{\widehat \sigma^2}{\widetilde \sigma_{1,1}^2} \biggr)^{1/\alpha^*}
\Biggr)
\Biggr\}
\\
&\qquad
+\frac{\sqrt{n}}{3\pi^2}
\bigl( \widehat \sigma^{2/\alpha^*} - (\sigma^*)^{2/\alpha^*} \bigr)
\bigl( \widetilde \sigma_{1,2}^{-2/\alpha^*} 
-\widetilde \sigma_{1,1}^{-2/\alpha^*} \bigr)
\\
&\qquad+ 
\frac{\sqrt{n}(\sigma^*)^{2/\alpha^*}}{3\pi^2}
\Bigl( 
\bigl( \widetilde \sigma_{1,2}^{-2/\alpha^*} 
-\widetilde \sigma_{1,1}^{-2/\alpha^*} \bigr)
- \bigl( (\sigma_{1,2}^*)^{-2/\alpha^*} - (\sigma_{1,1}^*)^{-2/\alpha^*} \bigr)
\Bigr)
\\
&=
\frac{\sqrt{n}(\sigma^*)^{2/\alpha^*}}{3\pi^2}
\Bigl( 
\bigl( \widetilde \sigma_{1,2}^{-2/\alpha^*} 
-\widetilde \sigma_{1,1}^{-2/\alpha^*} \bigr)
- \bigl( (\sigma_{1,2}^*)^{-2/\alpha^*} - (\sigma_{1,1}^*)^{-2/\alpha^*} \bigr)
\Bigr)
+\oo_\PP(1)
\end{align*}
in the same way as above. Let $F(x,y) = y^{-1/\alpha^*} - x^{-1/\alpha^*}$ for $x,y>0$. 
Since it follows that
\begin{align*}
&\sqrt{n}
\Bigl( 
\bigl( \widetilde \sigma_{1,2}^{-2/\alpha^*} 
-\widetilde \sigma_{1,1}^{-2/\alpha^*} \bigr)
- \bigl( (\sigma_{1,2}^*)^{-2/\alpha^*} - (\sigma_{1,1}^*)^{-2/\alpha^*} \bigr)
\Bigr)
\\
&=
\sqrt{n} \bigl( F(\widetilde \sigma_{1,1}^2, \widetilde \sigma_{1,2}^2) 
-F((\sigma_{1,1}^*)^2, (\sigma_{1,2}^*)^2) \bigr)
\\
&= 
\pd F( (\sigma_{1,1}^*)^2, (\sigma_{1,2}^*)^2 )
\sqrt{n} 
\begin{pmatrix}
\widetilde \sigma_{1,1}^2 - (\sigma_{1,1}^*)^2 
\\
\widetilde \sigma_{1,2}^2 - (\sigma_{1,2}^*)^2
\end{pmatrix}
+\oo_\PP(1)
\\
&=
\frac{1}{\alpha^*} 
\bigl((\sigma_{1,1}^*)^{-2/\alpha^*-2}, -(\sigma_{1,2}^*)^{-2/\alpha^*-2} \bigr)
\sqrt{n} 
\begin{pmatrix}
\widetilde \sigma_{1,1}^2 - (\sigma_{1,1}^*)^2
\\
\widetilde \sigma_{1,2}^2 - (\sigma_{1,2}^*)^2
\end{pmatrix}
+\oo_\PP(1)
\\
& = \frac{(\sigma^*)^{-2/\alpha^*-2}}{\alpha^*}
((\lambda_{1,1}^*)^{1+\alpha^*}, -(\lambda_{1,2}^*)^{1+\alpha^*} )
\sqrt{n} 
\begin{pmatrix}
\widetilde \sigma_{1,1}^2 - (\sigma_{1,1}^*)^2
\\
\widetilde \sigma_{1,2}^2 - (\sigma_{1,2}^*)^2
\end{pmatrix}
+\oo_\PP(1),
\end{align*}
we have
\begin{align*}
\sqrt{n} (\widetilde \theta_2 - \theta_2^*) 
= 
\frac{(\sigma^*)^{-2}}{3\pi^2 \alpha^*}
((\lambda_{1,1}^*)^{1+\alpha^*}, -(\lambda_{1,2}^*)^{1+\alpha^*} )
\sqrt{n} 
\begin{pmatrix}
\widetilde \sigma_{1,1}^2 - (\sigma_{1,1}^*)^2
\\
\widetilde \sigma_{1,2}^2 - (\sigma_{1,2}^*)^2
\end{pmatrix}
+\oo_\PP(1).
\end{align*}
We also obtain
\begin{align*}
\sqrt{n} (\widetilde \theta_1 - \theta_1^*) 
&= \sqrt{n}(\widehat \kappa - \kappa^*) \widetilde \theta_2
+\kappa^* \sqrt{n}(\widetilde \theta_2 - \theta_2^*)
\\
&= 
\frac{\kappa^* (\sigma^*)^{-2}}{3 \pi^2 \alpha^*}
((\lambda_{1,1}^*)^{1+\alpha^*}, -(\lambda_{1,2}^*)^{1+\alpha^*} )
\sqrt{n} 
\begin{pmatrix}
\widetilde \sigma_{1,1}^2 - (\sigma_{1,1}^*)^2
\\
\widetilde \sigma_{1,2}^2 - (\sigma_{1,2}^*)^2
\end{pmatrix}
+\oo_\PP(1),
\\
\sqrt{n} (\widetilde \eta_1 - \eta_1^*) 
&= \sqrt{n}(\widehat \eta - \eta^*) \widetilde \theta_2
+\eta^* \sqrt{n}(\widetilde \theta_2 - \theta_2^*)
\\
&= \frac{\eta^* (\sigma^*)^{-2}}{3 \pi^2 \alpha^*}
((\lambda_{1,1}^*)^{1+\alpha^*}, -(\lambda_{1,2}^*)^{1+\alpha^*} )
\sqrt{n} 
\begin{pmatrix}
\widetilde \sigma_{1,1}^2 - (\sigma_{1,1}^*)^2
\\
\widetilde \sigma_{1,2}^2 - (\sigma_{1,2}^*)^2
\end{pmatrix}
+\oo_\PP(1),
\\
\sqrt{n} (\widetilde \sigma^2 - (\sigma^*)^2) 
&= \sqrt{n} \biggl(\frac{\widehat \sigma^2}{\widehat \theta_2} 
- \frac{(\sigma^*)^2}{\theta_2^*} \biggr) \widetilde \theta_2
+\frac{(\sigma^*)^2}{\theta_2^*} \sqrt{n}(\widetilde \theta_2 - \theta_2^*)
\\
&= \frac{(\theta_2^*)^{-1}}{3 \pi^2 \alpha^*}
((\lambda_{1,1}^*)^{1+\alpha^*}, -(\lambda_{1,2}^*)^{1+\alpha^*} )
\sqrt{n} 
\begin{pmatrix}
\widetilde \sigma_{1,1}^2 - (\sigma_{1,1}^*)^2
\\
\widetilde \sigma_{1,2}^2 - (\sigma_{1,2}^*)^2
\end{pmatrix}
+\oo_\PP(1).
\end{align*}
Let $v_1 = (3\pi^2, \kappa^*, \eta^*, 1, (\sigma^*)^2/\theta_2^*)^\TT$
and $v_2 = (0, \kappa^*, \eta^*, 1, (\sigma^*)^2/\theta_2^*)^\TT$.
We then find from \eqref{eq-3-1} that
\begin{align*}
\sqrt{n}
\begin{pmatrix}
\widetilde \theta_0 - \theta_0^*
\\
\widetilde \theta_1 - \theta_1^*
\\
\widetilde \eta_1 - \eta_1^*
\\
\widetilde \theta_2 - \theta_2^*
\\
\widetilde \sigma^2 - (\sigma^*)^2
\end{pmatrix}
&=
\frac{(\sigma^*)^{-2}}{3 \pi^2 \alpha^*}
((\lambda_{1,1}^*)^{1+\alpha^*} v_1, -(\lambda_{1,2}^*)^{1+\alpha^*} v_2)
\sqrt{n} 
\begin{pmatrix}
\widetilde \sigma_{1,1}^2 - (\sigma_{1,1}^*)^2
\\
\widetilde \sigma_{1,2}^2 - (\sigma_{1,2}^*)^2
\end{pmatrix}
+\oo_\PP(1)
\\
&\dto N(0,\mathcal J),
\end{align*}
where $\mathcal J$ is represented by using 
$\theta_{-1}^* = (\theta_1^*, \eta_1^*, \theta_2^*, (\sigma^*)^2)$ as follows.
\begin{align*}
\mathcal J &= 
\frac{2(\sigma^*)^{-4}}{9 \pi^4 (\alpha^*)^2}
((\lambda_{1,1}^*)^{1+\alpha^*} v_1, -(\lambda_{1,2}^*)^{1+\alpha^*} v_2)
\begin{pmatrix}
(\sigma_{1,1}^*)^4 & 0
\\
0 & (\sigma_{1,2}^*)^4 
\end{pmatrix}
((\lambda_{1,1}^*)^{1+\alpha^*} v_1, -(\lambda_{1,2}^*)^{1+\alpha^*} v_2)^\TT
\\
&=
\frac{2}{9\pi^4 (\alpha^*)^2} 
((\lambda_{1,1}^*)^2 v_1 v_1^\TT + (\lambda_{1,2}^*)^2 v_2 v_2^\TT )
\\
&= 2
\begin{pmatrix}
\frac{(\lambda_{1,1}^*)^2}{(\alpha^*)^2} & 
\frac{(\lambda_{1,1}^*)^2}{3\pi^2\theta_2^* (\alpha^*)^2}\theta_{-1}^*
\\
\frac{(\lambda_{1,1}^*)^2}{3\pi^2\theta_2^* (\alpha^*)^2} (\theta_{-1}^*)^\TT &
\frac{(\lambda_{1,1}^*)^2 + (\lambda_{1,1}^*)^2}{9\pi^4 (\theta_2^*)^2 (\alpha^*)^2}
(\theta_{-1}^*)^\TT \theta_{-1}^*
\end{pmatrix}
\\
&= 2
\begin{pmatrix}
J_{1,1} & J_{1,2}
\\
J_{1,2}^\TT & J_{2,2}
\end{pmatrix}
.
\end{align*}

\textit{Proof of (i) of Theorem \ref{th4-2}. }
One can show that under 
$\frac{n}{(M_1 \land M_2)^{2(\alpha^* \tand 1)}} \to 0$, 
\begin{equation*}
\widetilde\sigma_{l_1,l_2}^2
-\sum_{i=1}^{n}
(x_{l_1,l_2}(\widetilde t_i)-x_{l_1,l_2}(\widetilde t_{i-1}))^2
=\oo_\PP(1),
\end{equation*}
in the same way as above, and thus one has $\widetilde \theta_0 \pto \theta_0^*$.

\subsection{Proofs of the results in Subsection \ref{sec4-2}}

Since the proofs of Theorems \ref{th4-3} and \ref{th4-4} 
are same as those of Theorems \ref{th4-1} and \ref{th4-2}, 
we only show that the asymptotic variance is given by $\mathcal K$.

We have
\begin{align*}
\sqrt{n} (\overline \mu_0 -\mu_0^*)
&= \sqrt{n} \biggl\{
\biggl(\frac{\widecheck \sigma^2}{\overline \varsigma_{1,1}^2}\biggr)^{1/\widehat \alpha}
- \biggl(\frac{(\sigma^*)^2}{(\varsigma_{1,1}^*)^2}\biggr)^{1/\alpha^*}
\biggr\}
\\
&= (\sigma^*)^{2/\alpha^*} 
\sqrt{n} \bigl( \overline \varsigma_{1,1}^{2/\alpha^*} 
- (\varsigma_{1,1}^*)^{2/\alpha^*} \bigr)
+\oo_\PP(1)
\\
&= -\frac{(\sigma^*)^{2/\alpha^*}}{\alpha^*}
(\varsigma_{1,1}^*)^{-2/\alpha^*-2} 
\sqrt{n} \bigl( \overline \varsigma_{1,1}^2 - (\varsigma_{1,1}^*)^2 \bigr)
+\oo_\PP(1)
\\
&= -\frac{1}{\alpha^* (\sigma^*)^2}((\mu_{1,1}^*)^{1+\alpha^*},0)
\sqrt{n} 
\begin{pmatrix}
\overline \varsigma_{1,1}^2 - (\varsigma_{1,1}^*)^2
\\
\overline \varsigma_{1,2}^2 - (\varsigma_{1,2}^*)^2
\end{pmatrix}
+\oo_\PP(1).
\end{align*}
Set $\Psi(x,y) = (3\pi^2 (y^{-1/\alpha^*} -x^{-1/\alpha^*})^{-1})^{\alpha^*}$.
Since
\begin{equation*}
\pd \Psi((\varsigma_{1,1}^*)^2,(\varsigma_{1,2}^*)^2)
= \frac{1}{3\pi^2} (-(\mu_{1,1}^*)^{1+\alpha^*}, (\mu_{1,2}^*)^{1+\alpha^*}),
\end{equation*}
we have
\begin{align*}
\sqrt{n} (\overline \theta_2 -\theta_2^*)
&= \frac{\theta_2^*}{(\sigma^*)^2}\sqrt{n} (\overline \sigma^2 -(\sigma^*)^2)
+\oo_\PP(1)
\\
&= \frac{\theta_2^*}{(\sigma^*)^2} 
\sqrt{n} \bigl( \Psi(\overline \varsigma_{1,1}^2,\overline \varsigma_{1,2}^2)
-\Psi((\varsigma_{1,1}^*)^2,(\varsigma_{1,2}^*)^2) \bigr)
+\oo_\PP(1)
\\
&= \frac{\theta_2^*}{(\sigma^*)^2} 
\pd \Psi((\varsigma_{1,1}^*)^2,(\varsigma_{1,2}^*)^2)
\sqrt{n} 
\begin{pmatrix}
\overline \varsigma_{1,1}^2 - (\varsigma_{1,1}^*)^2
\\
\overline \varsigma_{1,2}^2 - (\varsigma_{1,2}^*)^2
\end{pmatrix}
+\oo_\PP(1)
\\
&= \frac{\theta_2^*}{3 \pi^2(\sigma^*)^2} 
(-(\mu_{1,1}^*)^{1+\alpha^*}, (\mu_{1,2}^*)^{1+\alpha^*})
\sqrt{n} 
\begin{pmatrix}
\overline \varsigma_{1,1}^2 - (\varsigma_{1,1}^*)^2
\\
\overline \varsigma_{1,2}^2 - (\varsigma_{1,2}^*)^2
\end{pmatrix}
+\oo_\PP(1),
\\
\sqrt{n} (\overline \theta_1 -\theta_1^*)
&= \frac{\theta_1^*}{3 \pi^2(\sigma^*)^2} 
(-(\mu_{1,1}^*)^{1+\alpha^*}, (\mu_{1,2}^*)^{1+\alpha^*})
\sqrt{n} 
\begin{pmatrix}
\overline \varsigma_{1,1}^2 - (\varsigma_{1,1}^*)^2
\\
\overline \varsigma_{1,2}^2 - (\varsigma_{1,2}^*)^2
\end{pmatrix}
+\oo_\PP(1),
\\
\sqrt{n} (\overline \eta_1 -\eta_1^*)
&= \frac{\eta_1^*}{3 \pi^2(\sigma^*)^2} 
(-(\mu_{1,1}^*)^{1+\alpha^*}, (\mu_{1,2}^*)^{1+\alpha^*})
\sqrt{n} 
\begin{pmatrix}
\overline \varsigma_{1,1}^2 - (\varsigma_{1,1}^*)^2
\\
\overline \varsigma_{1,2}^2 - (\varsigma_{1,2}^*)^2
\end{pmatrix}
+\oo_\PP(1),
\\
\sqrt{n} (\overline \sigma^2 -(\sigma^*)^2)
&= \frac{1}{3 \pi^2} 
(-(\mu_{1,1}^*)^{1+\alpha^*}, (\mu_{1,2}^*)^{1+\alpha^*})
\sqrt{n} 
\begin{pmatrix}
\overline \varsigma_{1,1}^2 - (\varsigma_{1,1}^*)^2
\\
\overline \varsigma_{1,2}^2 - (\varsigma_{1,2}^*)^2
\end{pmatrix}
+\oo_\PP(1).
\end{align*}
Therefore, setting
$v_3 = (3\pi^2/\alpha^*, \theta_1^*, \eta_1^*, \theta_2^*, (\sigma^*)^2)^\TT$ and
$v_4 = (0, \theta_1^*, \eta_1^*, \theta_2^*, (\sigma^*)^2)^\TT$,
we see
\begin{align*}
\sqrt{n}
\begin{pmatrix}
\overline \mu_0 - \mu_0^*
\\
\overline \theta_1 - \theta_1^*
\\
\overline \eta_1 - \eta_1^*
\\
\overline \theta_2 - \theta_2^*
\\
\overline \sigma^2 - (\sigma^*)^2
\end{pmatrix}
&=
\frac{1}{3 \pi^2 (\sigma^*)^2}
(-(\mu_{1,1}^*)^{1+\alpha^*} v_3, (\mu_{1,2}^*)^{1+\alpha^*} v_4)
\sqrt{n} 
\begin{pmatrix}
\overline \varsigma_{1,1}^2 - (\varsigma_{1,1}^*)^2
\\
\overline \varsigma_{1,2}^2 - (\varsigma_{1,2}^*)^2
\end{pmatrix}
+\oo_\PP(1)
\\
&\dto N(0,\mathcal K),
\end{align*}
where $\mathcal K$ is expressed as follows.
\begin{align*}
\mathcal K &= 
\frac{2}{9 \pi^4 (\sigma^*)^4}
((\mu_{1,1}^*)^{1+\alpha^*} v_3, -(\mu_{1,2}^*)^{1+\alpha^*} v_4)
\begin{pmatrix}
(\varsigma_{1,1}^*)^4 & 0
\\
0 & (\varsigma_{1,2}^*)^4 
\end{pmatrix}
((\mu_{1,1}^*)^{1+\alpha^*} v_3, -(\mu_{1,2}^*)^{1+\alpha^*} v_4)^\TT
\\
&=
\frac{2}{9\pi^4} 
((\mu_{1,1}^*)^2 v_3 v_3^\TT + (\mu_{1,2}^*)^2 v_4 v_4^\TT )
\\
&= 2
\begin{pmatrix}
\frac{(\mu_{1,1}^*)^2}{(\alpha^*)^2} & 
\frac{(\mu_{1,1}^*)^2}{3\pi^2 \alpha^*}\theta_{-1}^*
\\
\frac{(\mu_{1,1}^*)^2}{3\pi^2 \alpha^*} (\theta_{-1}^*)^\TT &
\frac{(\mu_{1,1}^*)^2 + (\mu_{1,2}^*)^2}{9\pi^4}
(\theta_{-1}^*)^\TT \theta_{-1}^*
\end{pmatrix}
\\
&= 2
\begin{pmatrix}
K_{1,1} & K_{1,2}
\\
K_{1,2}^\TT & K_{2,2}
\end{pmatrix}
.
\end{align*}



\end{document}